\documentclass[a4paper,10pt]{amsart}

\usepackage{amssymb, amscd}
\usepackage{latexsym,epsfig}
\usepackage[all]{xy}
\usepackage{pst-all}
\usepackage{fancyhdr}
\numberwithin{equation}{section}
\def\today{\ifcase\month\or Jan\or Febr\or  Mar\or  Apr\or May\or Jun\or  Jul\or Aug\or  Sep\or  Oct\or Nov\or  Dec\or\fi \space\number\day, \number\year}

\newcommand\eprint[1]{Eprint:~\texttt{#1}}

\newcommand{\EE}{\mathbb E}
\newcommand{\FF}{\mathbb F}
\newcommand{\GG}{\mathbb G}
\newcommand{\HH}{\mathbb H}

\newcommand{\PP}{\mathbb P}
\newcommand{\QQ}{\mathbb Q}

\newcommand{\ZZ}{\mathbb Z}

\def\cH{{\mathcal H}}
\def\cF{{\mathcal F}}
\def\cK{{\mathcal K}}
\newcommand\Weyl[2]{{W}^{#1}_{#2}}
\newcommand\Weylp[2]{{W}'^{#1}_{#2}}

\newcommand\riso{\mathrel{\hskip2pt\raise-2.5pt\hbox{$\widetilde{\phantom{xx}}$}
\kern-16pt\longrightarrow}}
\newcommand\liso{\mathrel{\hskip2pt\raise-2.5pt\hbox{$\widetilde{\phantom{xx}}$}
\kern-17pt\longleftarrow}}
\newcommand\mapright[1]{\stackrel{#1}{\longrightarrow}}

\newcommand\rk{\mathop{\rm rk}\nolimits}

\newcommand\Legendre[2]{{\genfrac{(}{)}{}{}{#1}{#2}}}
\newcommand\Fl{\mathcal{F}\ell}

\newcommand\NS{\rm {NS}}
\newcommand\gr{\rm {gr}}
\newcommand\disc{\symb{disc}}
\newcommand\dual[1]{{#1}^*}

\def\Tensor{\otimes}
\def\Uo{{\mathcal U}}
\def\Uc{\overline{\mathcal U}}
\def\Vo{{\mathcal V}}
\def\Vc{\overline{\mathcal V}}
\def\cV{\overline{\mathcal V}}
\def\Sp{\rm Sp}
\def\SO{\rm SO}
\def\Ogrp{\rm O}
\def\sE{\mathcal E}
\def\sF{\mathcal F}
\def\sG{\mathcal G}
\def\sH{\mathcal H}
\def\sL{\mathcal L}
\def\cL{\mathcal L}
\def\sM{\mathcal M}
\def\sO{\mathcal O}
\def\sV{\mathcal V}
\def\sD{\mathcal D}
\newcommand\Zflag{\mathcal{ZF}}
\newcommand\Zframe{\mathcal{FF}}
\def\disc{\rm disc}
\def\iso{\cong}
\def\tensor{\otimes}

\numberwithin{equation}{section}
\newtheorem{theorem}{Theorem}[section]
\newtheorem{lemma}[theorem]{Lemma}
\newtheorem{proposition}[theorem]{Proposition}
\newtheorem{corollary}[theorem]{Corollary}

\newtheorem{definition-lemma}[theorem]{Definition-Lemma}

\theoremstyle{definition}
\newtheorem{definition}[theorem]{Definition}
\newtheorem{example}[theorem]{Example}

\theoremstyle{remark}
\newtheorem{remark}[theorem]{Remark}

\newtheorem{conventions}[theorem]{Conventions}

\begin{document}

\title[Cycle Classes on the Moduli of K3 surfaces in positive characteristic]{Cycle Classes on the Moduli of K3 surfaces \\ in positive characteristic}
\author{Torsten Ekedahl}
\thanks{Torsten Ekedahl unexpectedly passed away on November 23, 2011}
\address{Department of Mathematics\\ Stockholm University\\
 SE-106 91  Stockholm\\ Sweden}
\email{teke@math.su.se}

\author{Gerard van der Geer} \address{Korteweg-de Vries Instituut\\
Universiteit van Amsterdam\\ Postbus 94248\\1090 GE Amsterdam
The Netherlands}

\email{geer@science.uva.nl}


\begin{abstract}
This paper provides explicit closed formulas in terms of tautological classes
 for the cycle classes of the height and Artin invariant strata in families of
  K3 surfaces. The proof is uniform for all strata and uses a flag space as the
  computations in \cite{ekedahl10::cycle+class+e+o+strat} for the Ekedahl-Oort
  strata for families of abelian varieties, but employs a Pieri formula formula
  to determine the push down to the base space.
\end{abstract}

\maketitle

\begin{section}{Introduction}
Moduli spaces of algebraic varieties in positive characteristic possess
stratifications for which there are no analogues in characteristic zero.
This can make these moduli spaces more accessible than their counterparts in
characteristic $0$. 
The first example is the moduli space of elliptic curves where the distinction
ordinary versus supersingular provides a stratification. This
generalizes to the moduli of abelian varieties where one finds 
the Ekedahl-Oort stratification 
and the Newton polygon stratification. 
Besides the case of abelian varieties, 
where this phenomenon has attracted a lot of attention, the moduli of 
K3 surfaces in characteristic $p>0$ provide a beautiful example. 
To define a stratification one looks at typical characteristic $p$ invariants,
like the height and the Artin invariant.
The height appears if one considers multiplication by $p$ 
on the $1$-dimensional formal Brauer
group associated to the second \'etale cohomology group of a K3 surface; 
multiplication by $p$  is either zero or takes the form
$$
[p]\, t = a \, t^{p^h} +\text{higher order terms}
$$
with $a\neq 0$ and $t$ a local parameter. The number $h$ is called the
height and satisfies $1\leq h \leq 10$ or we have $[p]=0$ and then we put
$h=\infty$. If $h=\infty$ the
K3 surface is called supersingular (in the sense of Artin, see \cite{artin74::super+k3}). The loci of K3 surfaces with height 
$\geq h$ stratify the moduli. In the supersingular case one finds 
a further invariant, the Artin invariant $\sigma_0$, by looking at the
discriminant of the intersection pairing on the N\'eron-Severi group 
which is of the form $-p^{2\sigma_0}$ with $10 \geq \sigma_0 \geq 1$,
and this further stratifies the smallest stratum, the supersingular locus.
The generic supersingular case has $\sigma_0=10$ while $\sigma_0=1$ is
the most special case.
We thus find a stratification of $20$ strata on the 
$19$-dimensional moduli space, linearly ordered by inclusion. 

It is the purpose of this paper to calculate closed formulas for the cycle
classes of such strata in the Chow groups with rational coefficients
of moduli spaces of polarized K3 surfaces.  One may view such formulas
as a generalization of Deuring's formula that gives the number of
supersingular elliptic curves.
For the strata indexed by the height of the formal
Brauer group this was done in \cite{geer00::k3} in a somewhat ad hoc manner, but
for the more elusive strata parametrized by the Artin invariant this problem
remained open. It turns out that the cycle classes of the strata
can be expressed in powers of the tautological class $\lambda_1$, the first Chern
class of the Hodge bundle. The coefficients are complicated expressions
in $p$, the characteristic of our ground field. The result parallels
our joint work in \cite{ekedahl10::cycle+class+e+o+strat} that gives
such formulas for the moduli of abelian varieties.

The two invariants, the height and the Artin invariant, are of a seemingly different nature. However, they can be given a uniform description
in terms of the relative position of two filtrations
on the middle de Rham cohomology that refine the Hodge filtration and the
conjugate filtration; this is similar to how the Ekedahl-Oort strata 
on the moduli of abelian varieties, originally defined in terms
of group schemes, were interpreted 
in terms of relative position of flags on de Rham cohomology
in \cite{G99,ekedahl10::cycle+class+e+o+strat}. 
These two filtrations form a socalled $F$-zip in
the sense of \cite{moonen04::discr}.
The failure of transversality of these two filtrations
is measured by a double coset of a Weyl group and gives rise
to discrete invariants, like the height of the formal Brauer group and the
Artin invariant. In contrast, in characteristic $0$ the Hodge filtration 
and its complex conjugate are  always transversal. The role of the first 
cohomology group for abelian varieties is replaced by the second 
cohomology group. We shall explain the precise relation between
the discrete invariants and the relative positions of the two filtrations.

We consider here the moduli of lattice polarized K3 surfaces, i.e.,  
K3 surfaces together with an
embedding of a non-degenerate lattice in the N\'eron-Severi group
and then consider the primitive cohomology,  that is, the orthogonal complement of our non-degenerate lattice. 
The second de Rham cohomology of a K3 surface in characteristic $p$
comes with two filtrations, the Hodge filtration and the conjugate filtration.
We show that these two filtrations can be refined to a so-called pair of 
complete self-dual filtrations on the primitive cohomology
which are compatible with respect to the action of Frobenius; they form
 a so-called complete $F$-zip. 
Such a refinement is not unique, but the relative position is controlled by
an element of a Weyl group and this element is unique.

This naturally forces us to work in a flag space ${\mathcal F}$ 
of the primitive part of the second de Rham cohomology over a 
moduli space ${\mathcal M}$ of lattice polarized K3 surfaces. 
Looking at the full filtrations we find a stratification on 
${\mathcal F}$ indexed by elements in a Weyl group. The strata in 
${\mathcal F}$ project to strata on ${\mathcal M}$ and we are 
interested in formulas for the cycle classes of the strata on 
${\mathcal M}$. The reason that we nevertheless insist on working
on the flag space ${\mathcal F}$ is that the strata on ${\mathcal F}$
are much better behaved than on ${\mathcal M}$. Locally on
${\mathcal F}$ one can compare the strata with the Schubert strata
on the space of complete self-dual flags on an orthogonal space.
The idea, already used in the analogous situation for abelian varieties in \cite{ekedahl10::cycle+class+e+o+strat},
is that our flag space ${\mathcal F}$ as a stratified space
at a point can be identified up to the $(p-1)$st neighborhood with the flag space at an appropriate point.
This provides a lot of information about local structure of our strata
(dimension, Cohen-Macaulay-ness). It turns out that the strata
corresponding to special elements in the Weyl group (called final elements)
map in a finite surjective \'etale way to strata on ${\mathcal M}$.

As mentioned above, the strata on the moduli of polarized K3
surfaces (given by the height and Artin invariant) 
are linearly ordered. However, on the space of $1$-dimensional
isotropic subspaces  of a $n$-dimensional orthogonal 
vector space with $n$ even, 
the poset of Schubert varieties is not a total order; there are two 
mid-dimensional incomparable Schubert varieties (which are permuted under the orthogonal group). This points to a delicate subtlety here.
In the even-dimensional case one of the middle dimensional strata is excluded in the
strata on the moduli space ${\mathcal M}$. The reason behind this is
the existence of a deformation invariant, the Hodge
discriminant which, depending on its value (in ${\bf F}^*_p/{\bf F}^{*2}_p$),
excludes one of the mid-dimensional strata and leaves us with a
linearly ordered set of strata. 

This is reflected in the algebraic groups behind the scene.
In the case of abelian varieties it was shown \cite{ekedahl10::cycle+class+e+o+strat} that
for computing the classes the algebraic group ${\Sp}(2g)$ played an essential
r\^ole. By analogy with the complex case 
one would perhaps expect that in the case of
K3 surfaces the special orthogonal group ${\SO}(n)$ would play a similar r\^ole.
This is almost but not quite the case, it turns out that it is rather the full
orthogonal group ${\Ogrp}(n)$ that governs the situation. When the dimension $n$
(of the primitive part of cohomology) is odd the distinction between ${\SO}(n)$
and ${\Ogrp}(n)$ is not really seen (essentially as ${\Ogrp}(n)$ acts trivially on
the Dynkin diagram). The case of an even $n$ is markedly different (as this time
${\Ogrp}(n)$ acts non-trivially on the Dynkin diagram). 

Apart from these complications that appear when $n$ is even, our strategy for
finding cycle class formulas is the same whether $n$ is even or odd. Just as for
the case of abelian varieties 
(cf., \cite{ekedahl10::cycle+class+e+o+strat}) we work with a
space of complete flags extending the Hodge filtration and first obtain formulas
there for the classes of strata that are in bijection with the Schubert cells on
the complete flag space (of ${\SO}(n)$). We then push down these formulas to the
moduli space under consideration.

At this point however we follow a strategy which is different from that
used in \cite{ekedahl10::cycle+class+e+o+strat}: instead of using formulas of
Fulton and Pragacz we shall use a Pieri type formula. This fits well with the
fact that we have linearly ordered strata, but it introduces a new
problem. This Pieri formula involves many different strata all of which will
have to be pushed down to the moduli space. Comparing with the map from the
complete flag space $\Fl_n$ (of ${\SO}(n)$) to the space ${\mathcal I}(n)$ 
of isotropic $1$-dimensional subspaces we have the following situation. 
In the case of $\Fl_n$ each Schubert cell of the complete flag space maps 
to a Schubert cell of ${\mathcal I}(n)$. For each Schubert cell of 
${\mathcal I}(n)$ there is a unique Schubert cell
of $\Fl_n$ that maps isomorphically to it (the \emph{final} cell in our
terminology). All non-final cells map to a cell with positive dimensional fibres
and hence their cycle classes will push down to $0$. In our case the situation
is the same up to infinitesimal order $p-1$. That means that the map on a final
stratum is \'etale and on a non-final stratum it is non-separable. 
The degree with
which a final stratum maps to a stratum on our moduli space can be computed (and
usually is greater than $1$) and the result is analogous to the case of abelian varieties.
For a non-final stratum we can either see that its image is lower-dimensional,
and hence can be ignored, or we can find a factorization, called a
\emph{shuffle}, of the projection as an inseparable map of computable degree to
another stratum and the projection of that latter stratum.  Iterating this we
either get that a stratum has lower-dimensional image or that the projection
factors as an inseparable map (of computable degree) to a final stratum and the
projection of the final stratum. This allows us to get a complete description of
the push down of the classes coming from the Pieri formula and thus we get 
formulas for the cycle classes of the strata on the moduli space.

To give a feeling for the resulting formulas let us consider the simple case of
the moduli space ${\mathcal M}_d$ of K3 surfaces with a polarization 
of degree $d$, prime to the characteristic of the field $k$. 
One has $20$ strata $\Vo_w$ parametrized by
so-called final elements $w=w_i$ with $i=1,\ldots, 20$ in a Weyl group. The
strata $\Vo_{w_j}$ for $j=1,\ldots,10$ are the strata of K3 surfaces whose
formal Brauer group has finite height $j$, 
the stratum $\Vo_{w_{11}}$ is the
supersingular stratum and the strata indexed by $w_j$ for $j=12,\ldots,20$
correspond to supersingular 
K3 surfaces with Artin invariant $21-j$.  The strata 
come with a natural scheme structure. Our result expresses the
cycle classes of these strata as multiples of 
powers of the Hodge class  
$\lambda_1=c_1(\pi_* \Omega^2_{{\mathcal X}/M})\in 
{\rm CH}_{\QQ}^1({\mathcal M}_d)$
with  $\pi\colon {\mathcal X} \to {\mathcal M}_d$ the universal
K3 surface. 

\begin{theorem}\label{examplethm}
The cycle classes of the final strata $\Vc_w$ on the moduli space 
${\mathcal M}_d$ are
polynomials in $\lambda_1$ with coefficients that are $1/2$ times an integral
polynomial in $p\neq 2$ given by
\begin{eqnarray*}
{\rm i)} \quad 
[\Vc_{w_k}] &=& (p-1)(p^2-1)\cdots(p^{k-1}-1) \lambda_1^{k-1} \quad
\hbox{\rm if $1\leq k\leq 10$,}\\
{\rm ii)} \quad [\Vc_{w_{11}}] &=&\frac{1}{2} (p-1)(p^2-1)\cdots(p^{10}-1) 
\lambda_1^{10},\\
{\rm iii)} \quad 
[\Vc_{w_{10+k}}] &=&\frac{1}{2}
\frac{(p^{2k}-1)(p^{2(k+1)}-1)\cdots(p^{20}-1)}{(p+1)\cdots(p^{11-k}+1)}
\lambda_1^{9+k} \quad
\hbox{\rm if $2\leq k\leq 10$.}
\end{eqnarray*}
\end{theorem}
The appearance of the factor $1/2$ is related to the fact that the formulas of
  \cite[Thm 14.2 and Section 15]{geer00::k3} count the infinite height 
stratum doubly (cf.\ also \cite{geer01::formal+brauer}).
The moduli space is non-complete, but the formulas still make sense
on an appropriate compactification. 

Such formulas can be seen as a generalization of the well-known Deuring
formula for the number of isomorphism classes of 
supersingular elliptic curves over an algebraically closed field
of characteristic $p$.
The formulas for the height strata were already determined in \cite{geer00::k3}
in a completely different way, but that approach does not generalize
to the remaining strata.
The above theorem corresponds to the case where $n=2m+1$ is odd.
The more general case of moduli stacks of K3 surfaces with a marking of
a non-degenerate lattice in their N\'eron-Severi group forces us to treat
also the subtler case where $n$ is even. We finish this paper by giving two
examples that show that the even case appears quite naturally.

When dealing with K3 surfaces we do not have to go further than $n=21$.
However, our results should also be applicable to other moduli spaces related to
arithmetic subgroups of the orthogonal group, like moduli spaces of
hyperk\"ahler manifolds in positive characteristic which would give examples
with larger~$n$.

We expect that as the moduli of K3 surfaces in positive characteristic  gradually
become better understood our formulas will find many applications. Here
we give two
applications, one to the non-existence of supersingular 
elliptic K3 surfaces with a section with $\sigma_0=10$, see Proposition
\ref{sigma0neq10},
 and a
similar one to Enriques surfaces.

After the first version of this paper appeared important advances have been 
made concerning K3 surfaces. 
We mention the results of Maulik \cite{M12}, Madapusi Pera \cite{MP12,MP13} 
and Charles \cite{Ch13}
on Artin's conjecture and the Tate conjecture and results
of Liedtke \cite{L13} concerning the unirationality and 
moduli of supersingular K3 surfaces.
Artin's conjecture says that supersingular K3 surfaces 
(that is, of height $h=\infty$) have Picard number $\rho=22$ 
(that is, are supersingular in Shioda's sense). 
This has been verified for $p\geq 5$; 
the case $p=2$ was already done by Rudakov and Shafarevich, \cite{RS}. 
The papers of Madapusi Pera and Liedtke contain important information 
concerning moduli spaces of K3 surfaces in positive characteristic.  

\begin{conventions}
Throughout this paper we assume that the characteristic $p$ is not $2$
as orthogonal groups show a different behavior in characteristic $2$.
\end{conventions}

The original version of this paper was put on arXiv on April 15, 2011 
({\tt arXiv:1104.3024v1}).
On November 23, 2011 Torsten Ekedahl suddenly died. The second author has
revised this paper trying to improve the exposition and 
add more explanation. The mathematical content is essentially the same as in 
the original preprint. 

\end{section}
\begin{section}{Combinatorics}
We start with an auxiliary section on the combinatorics of the Weyl groups
associated to our orthogonal groups. We distinguish the B, C and  D cases.
As a general reference the reader might use \cite{bourbaki::Lie} or also
\cite{BB05}.
\begin{subsection}{B and C combinatorics}

  The Weyl group $\Weyl{B}{m}$ of ${\SO}(2m+1)$ can be identified with the
  subgroup of $S_{2m+1}$, the symmetric group on $2m+1$ letters, consisting of
  the permutations $\sigma \in S_{2m+1}$, for which
  $\sigma(i)+\sigma(2m+2-i)=2m+2$. We shall specify such a permutation by giving
  the images of the $1 \leq i \leq m$ as $[a_1,a_2,\ldots,a_m]$.  Thus the
  condition that this specify an element of $\Weyl{B}{m}$ is that $a_i \notin
  \{ a_j,m+1,2m+2-a_j\}$ for all $i \neq j$.  The elements which are reduced with
  respect to the set of roots obtained by removing the first root (so that the
  remaining roots form a root system of type $B_{m-1}$) are precisely those of
  the form $[a_1,a_2,\dots,a_m]$ with $a_1\neq m+1$ and $a_2,\ldots,a_m$ being
  an increasing sequence consisting of the first $m-1$ integers $\geq 1$ which
  are different from $a_1$ and $2m+2-a_1$ (cf.,
  \cite[\S3.4]{billey00::singul+schub}).  We write them as
  $[2m+2-a,1,2,3,\dots]$ including of course examples such as
  $[2,1,3,\dots]$ and $[2m+1,2,3,\dots]$. We shall call these elements the
  {\sl final elements} of $\Weyl{B}{m}$. There are $2m$ final elements;
we shall list these looking for the next largest $w(1)$ and thus
denote these by $w_1=[2m+1,2,3,\dots]$, $w_2=[2m,1,3,\ldots]$ (if $m\neq 1$),
$\ldots,w_{2m}=1$ and we sometimes write $w_{\emptyset}$ for $w_1$.

The {\sl simple reflections} $s_i$ for $i=1,\dots,m$ of $\Weyl{B}{m}$ are
the permutations $s_i=(i,i+1)(2m+1-i,2m+2-i)$ for $i=1,\dots,m-1$ and
$s_m=(m,m+2)$. We also define the {\sl weight representation} of
$\Weyl{B}{m}$ on ${\ZZ}^m$ with basis vectors $\epsilon_i$ ($i=1,\ldots,m$) 
given by, for $\sigma \in \Weyl{B}{m}$,
\begin{displaymath}
\sigma(\epsilon_i)=
\begin{cases}
\epsilon_{\sigma(i)} & \text{ if }\sigma(i)\leq m \text{ and}\\
-\epsilon_{2m+2-\sigma(i)} & \text{ if }\sigma(i) > m.
\end{cases}
\end{displaymath}
We thus can view $\Weyl{B}{m}$ as a reflection group of this lattice.  In
particular, for an element $\alpha \in {\ZZ}^m$ we have the reflection 
$s_{\alpha}$
in $\alpha$ with $s_{\alpha}(x)=x-\langle \alpha, x\rangle \alpha$; e.g.\
$s_i=s_{\epsilon_i}$.

For a permutation $w$ of $\{1,2,\dots,n\}$ we define
\begin{displaymath}
r_w(i,j)=\# \{ 1 \leq a \leq i : w(a) \leq j\}
\end{displaymath}
for $1\leq i,j\leq n$. It is clear that a permutation is determined by this
function.

The length of an element of $\Weyl{B}{m}$ (in the sense of Coxeter groups) may be
described in concrete terms as
\begin{displaymath}
\ell(w) = \# \{ 1\leq i\leq j\leq m: w(i)>w(j)\}+
\#\{ 1\leq i\leq j\leq m : w(i)+w(j)>2m+2 \}.
\end{displaymath}

We shall occasionally have to deal with the Weyl group $\Weyl{C}{m}$ of
${\Sp}(2m)$. It has almost exactly the same description as $\Weyl{B}{m}$,
in fact it is isomorphic to it, except
that it is considered as subgroup of $S_{2m}$:
$$
\Weyl{C}{m}:=\{ \sigma \in S_{2m} :  \sigma(i)+\sigma(2m+1-i)=2m+1\}\, ;
$$
a correspondence between them is given
by $w \in \Weyl{B}{m}$ defining an element $w' \in \Weyl{C}{m}$ by $w':= \sigma
w\sigma^{-1}$ where $\sigma(i)=i$ if $1 \leq i \leq m$ and $\sigma(i)=i-1$ if $m+1
< i \leq 2m+1$. The length of an element is given by
\begin{displaymath}
\ell(w) = \# \{1\leq i<j\leq m :  w(i)>w(j)\} +
\# \{1\leq i\leq j\leq m : w(i)+w(j)>2m+1\}.
\end{displaymath}
Finally, we define the {\sl discriminant}, ${\disc}(w)\in \{+1,-1\}$,
of $w \in \Weyl{B}{m}$ to be the sign of $w$ as an element of $S_{2m+1}$.
The reason for calling this homomorphism `${\disc}$' will appear later.
\end{subsection}

\begin{subsection}{D combinatorics}

  The Weyl group $\Weyl{D}{m}$ of ${\SO}(2m)$ consists of the permutations in
  $\sigma \in S_{2m}$ for which $\sigma(i)+\sigma(2m+1-i)=2m+1$ and such that
  there is an even number of $1\leq i \leq m$ for which $\sigma(i) > m$. The
  subgroup of $S_{2m}$ fulfilling the same conditions except for the parity
  condition form a subgroup of $S_{2m}$ which can be identified with the Weyl
  group $\Weyl{C}{m}$ for ${\Sp}(2m)$. Hence $\Weyl{D}{m}$ is a subgroup of
  $\Weyl{C}{m}$ of index $2$ and more precisely it is the kernel of the
  signature homomorphism ${\rm sign}: {\Weyl{C}{m}}\to {\pm 1}$.  We denote a
  permutation in $\Weyl{C}{m}$ as $[a_1,a_2,\dots,a_m]$. Thus the condition
  that this specify an element of $\Weyl{C}{m}$ is that $a_i \notin \{
  a_j,2n+1-a_j\}$ for all $i \neq j$ and it belongs to $\Weyl{D}{m}$ if also the
  number of $a_i$ with $a_i>m+1$ is even. The length of an element fulfils the
  formula
\begin{displaymath}
\ell(w) = \# \{1\leq i<j\leq m: w(i)>w(j)\}+
\# \{1\leq i < j\leq m :w(i)+w(j)>2m+1\}.
\end{displaymath}
The {\sl simple reflections} $s_i$, $i=1,\dots,m$, of $\Weyl{D}{m}$ are the
permutations $s_i=(i,i+1)(2m-i,2m+1-i)$ for $i=1,\dots,m-1$ and
$s_m=(m-1,m+1)(m,m+2)$. The simple reflections of $\Weyl{C}{m}$ are the $s_i$,
$i=1,\dots,m-1$, and $s'_m=(m,m+1)$. We also have the {\sl weight
representation} of $\Weyl{C}{m}$ on ${\ZZ}^m$ with basis vectors $\epsilon_i$
($i=1,\ldots,m$)  given by
\begin{displaymath}
\sigma(\epsilon_i)=
\begin{cases}
\epsilon_{\sigma(i)}&\text{ if } \sigma(i)\leq m \text{ and}\\
-\epsilon_{2m+1-\sigma(i)}&\text{ if }\sigma(i) >m. 
\end{cases}
\end{displaymath}
Note that the fact that the larger group is equal to $\Weyl{C}{m}$ is somewhat
accidental. To us it will rather be the Weyl group of ${\Ogrp}(2m)$ 
(as opposed to the Weyl group of ${\SO}(2m)$) or, equivalently, 
as the group generated by
$\Weyl{D}{m}$ and the non-trivial graph automorphism of $D_m$ (which in the
$D_4$ case is the one permuting the last two vertices). From the latter point of
view $s'_m$ gives a non-trivial graph automorphism, indeed it commutes with
$s_i$, $1\leq i <m-1$ and conjugation by it permutes $s_{m-1}$ and $s_m$. To
emphasize this point of view we shall, when relevant, write the supergroup
$\Weyl{C}{m}$ as $\Weylp{D}{m}$. In this context we need a definition of length
on $\Weylp{D}{m}$ that mimics the length of $\Weyl{D}{m}$ (rather than that of
$\Weyl{C}{m}$) and which we shall therefore denote $\ell_D$:
\begin{displaymath}
\ell_D(w) = \#\{1\leq i<j\leq m : w(i)>w(j)\}+
\# \{1\leq i < j\leq m : w(i)+w(j)>2m+1\}
\end{displaymath}
It has the property that its restriction to $\Weyl{D}{m}$ equals its natural length
and that $\ell_D(ws'_m)=\ell_D(w)$.

The elements which are reduced with respect to the set of roots of $D_m$
obtained by removing the first root (so that the remaining roots form a root
system of type $D_{m-1}$) are precisely those of the form $[a_1,a_2,\dots,a_m]$
with $(a_2,\ldots,a_m)$ being the lexicographically smallest sequence of
integers making $[a_1,a_2,\dots,a_m]$ an element of $\Weyl{D}{m}$. 
We list these by looking for the next largest $w(1)$ 
and thus have $w_1:=[2m,2,3,\dots,m-1,m+1]$, 
encountering of course examples such as
$[2,1,3,\dots,m]$ and $[m+1,1,2,\dots,m-1,m+2]$. We
shall call these elements the {\sl final elements} of
$\Weyl{D}{m}$. There are $2m$ final elements. The longest one is $w_1$ which we
shall also denote $w_\emptyset$. It has the reduced expression $s_1s_2\cdots
s_{m-2}s_{m-1}s_ms_{m-2}\cdots s_1$. We also put $w'_k := w_ks'_m$, and
we shall call these the
{\sl twisted final elements} with the alternative notation $w'_\emptyset$
for $w'_1$.
\end{subsection}
\end{section}

\begin{section}{Flag spaces}
Central in this paper are filtrations on the second the Rham 
cohomology of a K3 surface or on a primitive part of that. 
The intersection form makes this cohomology space
into a quadratic space. In this auxiliary section
we shall be interested in flags in a finite-dimensional orthogonal or symplectic
space and we start by recalling some well-known facts. 
We refer to \cite{billey00::singul+schub} or \cite{bourbaki::Lie}.

Let thus $V$ be an
$n$-dimensional vector space over a field ${\bf k}$ provided with a 
non-degenerate quadratic or symplectic form $\langle -, - \rangle $.  A flag 
$(0)=V_0\subset V_1\subset V_2 \subset \ldots \subset V_r$ of 
subspaces of $V$ is called {\sl isotropic} if the
restriction of the form to $V_r$ is zero.  We say that the flag is
{\sl maximal} if $r=k:=[n/2]$ (and hence $\dim(V_i)=i$). We can extend a
maximal flag to a {\sl self-dual complete flag} by putting
$V_{j}=V_{n-j}^{\perp}$.
The group ${\SO}(n)$ does not always acts transitively on complete flags.  
Two flags $V_\bullet$ and $V'_\bullet$ are in the same orbit under 
conjugation by ${\SO}(n)$ precisely when $\dim(V_k\cap V'_k) \equiv k \bmod 2$.

Now, given a complete flag $V_\bullet$ we may construct another complete flag
$V'_\bullet$ as follows: We let $V'_i=V_i$ for $i\neq k$ and then let $V'_k$ be
the unique maximal totally isotropic subspace containing $V_{k-1}$ and being
contained in $V_{k+1}$ that is distinct from $V_k$. As $V_k\cap V'_k=V_{k-1}$ we
see that $V_\bullet$ and $V'_\bullet$ are not conjugate under ${\SO}(n)$. 
We shall call $V'_\bullet$ the {\sl flip} of $V_\bullet$.

When the space is symplectic or $n$ is odd complete flags correspond precisely
to Borel subgroups of the symplectic or special orthogonal group; one associates
to a flag its stabiliser. The orthogonal even case is different however. To us
the main difference is the fact that ${\SO}(n)$ does not act transitively on
complete flags.

This leads us to introduce the notion of
{\sl self-dual almost complete flag} (when $n=2k$) which is specified by
an isotropic flag 
$(0)=V_0\subset V_1\subset V_2 \subset \ldots \subset V_{k-1}$
where $\dim V_i = i$ (and hence extended to a larger flag by putting
$V_{j}=V_{n-j}^{\perp}$ for $k+1\leq j\leq n$). 
If we let $\sF_n$ be the space of almost complete flags and $\sF'_n$ 
the space of complete flags, then the forgetful map $\sF'_n \to \sF_n$ is 
an \'etale double cover whose associated involution map $\sF'_n \to \sF'_n$ 
is given by the flip. Furthermore, ${\SO}(n)$
acts transitively on $\sF_n$ with stabilisers the Borel group of it. On the
other hand ${\Ogrp}(n)$ acts transitively both on $\sF_n$ and $\sF'_n$. The
stabilisers for the action on $\sF_n$ are subgroups of ${\Ogrp}(n)$ whose
intersection with ${\SO}(n)$ are Borel subgroups and which map surjectively onto
${\Ogrp}(n)/{\SO}(n)$ whereas the stabilisers on $\sF'_n$ are the
Borel subgroups of ${\SO}(n)$.

More generally if we have an orthogonal vector bundle $\sE \to X$ of constant
rank $n=2k$, then we have the bundle of almost complete flags $\sF(\sE) \to X$
and complete flags $\sF'(\sE) \to X$ and an \'etale double cover $\sF'(\sE) \to
\sF(\sE)$. This double cover is actually the pullback of a double cover of $X$:
it can be obtained by considering the quadric $Y$ in the ${\PP}^1$-bundle 
${\PP}(V_{k+1}/V_{k-1})$ defined by the orthogonal form;
it defines a double cover, the
{\sl discriminant double cover}, $\pi: Y \to X$.
Thus we get a morphism $\sF'(\sE) \to Y$ which fits into a cartesian diagram
\begin{displaymath}
\begin{CD}
\sF'(\sE) @>>> Y\\
@VVV @VVV\\
\sF(\sE) @>>> X.
\end{CD}
\end{displaymath}
The special properties of the even orthogonal case is the reason for the
relevance of the group $\Weylp{D}{m}$ as the following proposition shows.
It gives representatives for the orbits of pairs of flags.

It is well-known that the relative position of two flags can be measured
by an element of a Weyl group. We spell out the result for the cases that
interest us.

\begin{proposition}\label{Flags classification}
\item{1)} Let $e_1,\dots,e_{2m}$ be the standard basis of a symplectic space with
$\langle e_i, e_j \rangle =\delta_{i,2m+1-j}$ for $j\geq i$. 
The orbits of the action of ${\Sp}(2m)$ on pairs of totally isotropic 
complete flags in $2m$-dimensional space are in bijection with the 
elements of the Weyl group $\Weyl{C}{m}$. The element $w \in \Weyl{C}{m}$ 
corresponds to the orbit of $((\sum_{j\leq i}{\bf k}
e_{j}),(\sum_{j\leq i}{\bf k} e_{w^{-1}(j)}))$.

\item{2)}
Let $e_1,\dots,e_{2m+1}$ be the standard basis of an orthogonal space with
$\langle e_i, e_j \rangle =\delta_{i,2m+2-j}$. The orbits of the action of 
${\SO}(2m+1)$ on pairs of totally isotropic complete flags in 
$2m+1$-dimensional space are in
bijection with the elements of the Weyl group $\Weyl{B}{m}$. The element 
$w \in \Weyl{B}{m}$ corresponds to the orbit of 
$((\sum_{j\leq i}{\bf k} e_{j}),(\sum_{j\leq  i}{\bf k} e_{w^{-1}(j)}))$.
\item{3)}
Let $e_1,\dots,e_{2m}$ be the standard basis of an orthogonal space with
$\langle e_i, e_j \rangle = \delta_{i,2m+1-j}$.  
The orbits of the action of ${\SO}(2m)$ on
pairs of totally isotropic complete flags in $2m$-dimensional space are in
bijection with the elements of the group $\Weylp{D}{m}$. An element $w$ in
$\Weylp{D}{m}$ corresponds to the orbit of 
$((\sum_{j\leq i}{\bf k} e_{j}),(\sum_{j\leq i}{\bf k} e_{w^{-1}(j)}))$.  
If $(F_\bullet,D_\bullet)$ lies in the orbit
corresponding to $w$, then ${\disc}(w)=(-1)^d$, where $d=\dim(E_m\cap
D_m)$. Flipping the first flag changes the type from $w$ to $ws'_m$ 
and flipping the second changes it from $w$ to $s'_mw$.
\begin{proof}
The first and second part is of course well known, the third part perhaps less
so but in any case is just as easy to prove.
\end{proof}
\end{proposition}

We shall say that a basis such as in the proposition is {\sl adapted} to
the two flags. We shall also say that two complete flags are in
{\sl relative position} $w$ for $w$ in $\Weyl{C}{m}$, $\Weyl{B}{m}$ or
$\Weylp{D}{m}$ respectively if they belong to the orbit above associated to
$w$. Note that in the $B$ and $C$ cases we are dealing with orbits of $G$ 
(which equals ${\SO}(2m+1)$, resp.\ ${\Sp}(2m)$) on the product of 
flag spaces $G/B\times G/B$ and we are dealing with the well-known 
bijection between such orbits. In the even orthogonal case 
(and when ${\bf k}=\overline{\bf k}$) flags are in bijection
with ${\Ogrp}(2m)/B$, where $B$ is a Borel group of ${\SO}(2m)$ 
(the stabiliser of a fixed flag) where of course ${\Ogrp(}2m)/B$ 
has two components. Orbits under ${\Ogrp}(2m)$ of pairs of flags 
are then in bijection with $\Weylp{D}{m}$. We may however reduce 
ourselves at will to just the action of ${\SO}(2m)$ on
${\SO}(2m)/B$. Indeed, $V_\bullet$ and $U_\bullet$ are in relative 
position $w$ precisely when $V_\bullet$ and $U'_\bullet$ are in 
relative position $ws'_m$ where $U'_\bullet$ is the flip of $U_\bullet$.

All this relativizes to the situation of a symplectic or orthogonal vector
bundle $\sV$ of rank $n$ over a base $S$ (in which $2$ is
invertible). We can then construct the {\sl flag bundle} $\Fl(\sV)$ of
complete self-dual flags in $\sV$. In the even orthogonal this factors as above
through the {\sl discriminant cover} $\sD_{\sV}$. 
The involution associated to the discriminant cover 
extends to an involution of $\Fl(\sV)$ taking a flag to its flip. The
same terminology will be used for partial flags that contain a middle
dimensional member.

For later purposes the pointwise definition of flags to be in 
relative position $w$ does not suffice. We recall from 
\cite{ekedahl10::cycle+class+e+o+strat} the scheme theoretic definition: 
if we have two flags
over an affine scheme $Y$ we have two sections $s,t : Y \to T$,
where $T$ is a $G/B$-bundle with structure group $G$ 
for a semi-simple group $G$ and a Borel group $B$. 
Then for any element $w$ of the Weyl group of $G$ we define a 
(locally) closed subscheme ${\mathcal U}_w$ (resp.~$\overline{\mathcal U}_w$) 
of $Y$ in the following way. We choose locally 
(possibly in the \'etale topology) a trivialization of $T$ 
for which $t$ is a constant section. Then $s$ corresponds
to a map $Y \to G/B$ and we let ${\mathcal U}_w$ (resp.~$\overline{\mathcal U}_w$) 
be the inverse image of the $B$-orbit $BwB$ (resp.~of its closure in $G/B$).
Another trivialization will differ
by a map $Y \to B$;  as $BwB$ and its closure are $B$-invariant these
definitions are independent of the chosen trivializations
 and hence give global subschemes on $Y$. If $s$ and
$t$ have the property that $Y={\mathcal U}_w$, 
then we shall say that $s$ and $t$ are in
{\sl relative position $w$} and if $Y=\overline{\mathcal U}_w$ 
we shall say that $s$ and $t$ are in {\sl relative position $\leq w$}.

\end{section}
\begin{section}{$F$-zips}
The Hodge filtration and the conjugate filtration on the second de Rham cohomology
(or a primitive part of that) of a K3 surface form a so-called orthogonal $F$-zip. In this auxiliary
section we introduce the stack of flagged $F$-zips and certain substacks of it.

Recall (cf., \cite{moonen04::discr}) that an orthogonal or symplectic $F$-zip is a
tuple $(M,C^\bullet,D_\bullet,\varphi_\bullet)$, where $M$ is an orthogonal or
symplectic vector bundle over a base of positive characteristic, $0=C^r\subseteq
C^{r-1}\subseteq\cdots\subseteq C^0=M$ and $0=D_0\subseteq D_{1}\subseteq
\cdots\subseteq D_r=M$ are self-dual (not necessarily complete) flags on $M$ and
$\varphi_\bullet$ a collection of isomorphisms
$\varphi_i : F^*(C^{i}/C^{i+1}) \to D_{i+1}/D_{i}$ compatible with the
isomorphisms $C^{i+1}/C^{i}\riso (C^{r-i}/C^{r-i-1})^*$ and $D_{i+1}/D_{i}\riso
(D_{r-i}/D_{r-i-1})^*$ induced by the pairing. If the rank of $D_i$ has the
constant value $n_i$ we say that the $F$-zip is of {\sl dimension type} 
$\underline n=(n_r,n_{r-1},\dots,n_0)$. 
A {\sl flagged $F$-zip} is an $F$-zip together
with a complete self-dual flag $0=E_0\subset E_1\subset\cdots\subset E_n=M$
with $C^i=E_{n_i}$ where $n_i:=\rk(C^i)$ (and $\rk(E_i)=i$). We can use
$\varphi_\bullet$ to extend the $D$ flag to a complete self-dual flag
$G_\bullet$ by the condition that $D_{i+1}\subseteq G_{j'}\subseteq D_i$
when $C^i\subseteq E_j\subseteq C^{i+1}$, where $j'-j=n_{m-i-1}-n_i$, $m$
being the rank of $M$, and $G_{j'}/D_{i+1}=\varphi_i(F^*(E_j/C^i))$.

We now want to introduce the stack of flagged $F$-zips.  Starting from the
algebraic stacks ${\rm BO}(m)$ and ${\rm BSp}(m)$ of orthogonal resp.\ symplectic
vector bundles of rank $m$ (over ${\bf Z}/p$) one builds the algebraic stack $\Zflag$
of flagged $F$-zips (with $\rk D^i=n_i$ and just as we have fixed these ranks we
fix whether or not we have a symplectic or orthogonal bundle). If we only have
an incomplete (but still self-dual) flag extending $C^\bullet$ we shall speak of
a {\sl partially flagged $F$-zip} and we can do the same construction
getting a partial flag extending $D_\bullet$. A partially flagged $F$-zip is
{\sl stable} if for every $i$ and every $k$ we have that $D_i\cap
C^k+C^{k+1}$ is equal to some $C^j$ (or in relevant cases the middle part of
 the flip of $C^\bullet$), where $i$ and $k$ are chosen so that $D_i$, $C^k$ and
$C^{k+1}$ are defined; here relevant means whenever the dimension is even
and the partial flag involves the middle part. 
The canonical map of $\Zflag$ to the stack of orthogonal 
$F$-zips ${\mathcal Z}$ is relatively representable, so $\Zflag$ is an algebraic stack.

If the rank of the flagged orthogonal $F$-zip $(M,E,G)$ is even we can replace
both $E$ and $G$ by their flips and it is easy to see that we get a new
flagged $F$-zip which will be called the {\sl flip} of the flagged $F$-zip.

Fixing the rank, $n$, of $M$ and the flavor (symplectic or orthogonal) the
relative position of the flags $C^\bullet$ and $D_{\bullet}$ is, by
Proposition \ref{Flags classification}, described by an element $w$ of
$\Weyl{C}{n/2}$, $\Weyl{B}{(n-1)/2}$, and $\Weylp{D}{n}$ when the $F$-zip is
symplectic, orthogonal with $n$ odd, and orthogonal with $n$ even
respectively. We call the $F$-zip of {\sl type} $w$ in this case.

This defines locally closed substacks $\Zflag_w$ of $\Zflag$
consisting of those flags (of fixed flavor 
and ranks $\underline{n}$) of relative
position $w$, i.e., type $w$. We now also fix the sequence
$\underline{n}=(0=n_{0}<n_1<\cdots<n_r=n)$ where we demand that $\rk(D_i)=n_i$
(in particular self-duality forces $n_{i+1}-n_{i}=n_{r-i}-n_{r-i-1}$) and let
$w_\emptyset$ be the element of the appropriate Weyl group (as subgroup of
$S_n$) that takes the first $n_0$ integers to the last $n_0$ (in order), the
next $n_1-n_0$ integers to the $n_1-n_0$ last (still in order) and so on (that
is, it sends the interval $[n_i+1,n_{i+1}]$ to $[n+1-n_{i+1},n-n_i]$ preserving
order).

It will be useful to have an explicit scheme with a flagged $F$-zip of type $w$ over it 
that lies faithfully flat over $\Zflag_w$.

\textbf{Construction}: Given an element $w$ in the appropriate group we define a
flagged $F$-zip over the affine scheme 
${\rm Spec}({\bf F}_p[x_{ij}]_{1\leq i<j\leq n}/I)$ as follows, where
generators of the ideal $I$ are specified below:
\begin{itemize}
\item $e_i$, $i=1,\dots,n$, is a basis for $M$ with $\langle e_i, e_j \rangle=\delta_{i,n+1-j}$
for $i\leq j$.

\item $C^i$ has $e_1,\dots,e_i$ as a basis and $D_i$ has
$e_{w^{-1}(1)},\dots,e_{w^{-1}(i)}$ as a basis.

\item For $n_k<i\leq n_{k+1}$ we have
$\varphi_k(e_i)=e_{w^{-1}w_\emptyset(i)}+\sum_{w_\emptyset w(j)<i}x_{ij}e_{j}
\bmod D_{n_{r-k-1}}$.

\item The matrix $\mathrm{Id}_n+(x_{ij})$ is symplectic or orthogonal respectively with
respect to the scalar product of the basis $e_1,\dots,e_n$.

\item $x_{ij}=0$ unless $i \prec j$, where $i \prec j$ precisely when
$w^{-1}(i)>w^{-1}(j)$ for $n_\ell<j<i\leq n_{\ell+1}$ for some $\ell$ with
$n_\ell < n/2$.
\end{itemize}

When $n$ is odd, then we can define another flagged $F$-zip over ${\rm Spec}
({\bf F}_p[x_{ij}]_{1\leq i<j\leq n}/I)$ with the same definition except that we let
$$
\varphi(e_{(n+1)/2})=-e_{(n+1)/2}+\sum_{w_\emptyset
  w(j)<(n+1)/2}x_{(n+1)/2,j}e_{j}
$$ 
instead of $e_{(n+1)/2}+\sum_{w_\emptyset
  w(j)<(n+1)/2}x_{(n+1)/2,j}e_{j}$. Therefore, we define $Y_w$ as
${\rm Spec}({\bf F}_p[x_{ij}]_{1\leq i<j\leq n}/I)$ if $n$ is even and the disjoint union of
two copies of it when $n$ is odd. In both cases there is a flagged $F$-zip $\cF_n$
over $Y_w$. When $n$ is even it is the one constructed above. 
When $n$ is odd we have
the one with $\varphi(e_{(n+1)/2})=e_{(n+1)/2}+\cdots$ on one copy and the one
with $\varphi(e_{(n+1)/2})=-e_{(n+1)/2}+\cdots$ on the other (the flip). 
By construction
the two flags are everywhere of type $w$. This gives us a map $Y_w \to
\Zflag_w$.

\begin{proposition}\label{Versal $F$-zips}
The map $Y_w \to \Zflag_w$ is faithfully flat.
\begin{proof}
  By assumption there is a frame space $\Zframe_w \to \Zflag_w$ of bases of a
  versal flagged $F$-zip on $\Zflag_w$ adapted to the two flags. It is a group
  torsor and since the group scheme is flat also faithfully flat 
so that it is enough to show that the induced
  map $Y_w \to \Zframe_w$ is faithfully flat. There are
  functions $x_{ij}$ for $n_\ell<j\leq i\leq n_{\ell+1}$ such that
  $\varphi_\ell(e_i)=x_{ii}e_{w^{-1}w_\emptyset(i)}+\sum_{w_\emptyset
    w(j)<i}x_{ij}e_{j} \bmod D_{\ell}$, where $x_{ii}\neq 0$. Let $T$ consist of
  the diagonal automorphisms $e_i \mapsto t_ie_i$, where $t_i\cdot t_{n+1-i}=1$
  and $t_{(n+1)/2}=1$ if $n$ is odd. It transforms a $\varphi$ into another $F$-zip
  with $x_{ii}=t_{w^{-1}w_\emptyset(i)}^{-1}t^{p}_i$. As the endomorphism of $T$
  given by $(t_i) \mapsto (t_{w^{-1}w_\emptyset(i)}^{-1}t^{p}_i)$ is separable
  (inducing $-w^{-1}w_\emptyset(i)$ on the Lie algebra) we get that the map from
  the substack of $\Zframe_w$ with $x_{ii}=1$ for $i\neq(n+1)/2$ and
  $x_{ii}=\pm1$ if $i=(n+1)/2$ to $\Zframe_w$ is an equivalence and hence we may
  restrict to it.

It remains to show that we may remove the $x_{ij}$ with $n_\ell<j<i\leq
n_{\ell+1}$ and $w^{-1}(j)<w^{-1}(i)$. Under those assumptions, the change
of basis
$e'_i=e_i+\lambda e_j$, $e'_{\overline{\jmath}}=e_{\overline{\jmath}}-\lambda
e_{\overline{\imath}}$, with $\overline{x}=n+1-x$, preserves both flags. We now
have
\begin{displaymath}
  \varphi_k(e'_i)=\varphi_k(e_i+\lambda
  e_j)=e_{w^{-1}w_\emptyset(i)}+\sum_{w_\emptyset w(k)<i}x_{ik}e_{k}+
  \lambda^p \bigr(e_{w^{-1}w_\emptyset(j)}+\sum_{w_\emptyset w(\ell)<j}x_{i\ell}e_{\ell}\bigl)
\end{displaymath}
and we try to choose $\lambda$ such that the coefficient in front of
$e_{w^{-1}w_\emptyset(j)}$ is equal to zero. This gives a monic equation in
$\lambda$ with $\lambda^p$ as top term and hence defines a surjective finite flat
covering. We can repeat this construction in a way so that we take the largest
$i$ and $j$ first in order for subsequent operations not to reintroduce non-zero
coefficients in positions where they have been removed. At the end we get the
chosen $F$-zip on
$Y_w$ which shows fully faithful flatness as each step is fully faithfully flat.
\end{proof}
\end{proposition}

\end{section}
\begin{section}{The Hodge discriminant}
As alluded to in the Introduction there is a subtle invariant of the cohomology that prevents one of the
middle dimensional strata to turn up in the stratification of our moduli spaces. In this section we study
this invariant.

We begin now by introducing a discriminant which is a Hodge theoretic description of
Ogus' crystalline discriminant (defined under slightly more general
circumstances). For that we need the determinant of a complex
in the sense of Mumford and Knudsen (cf., \cite{knudsen76}). Recall that in
order to get the signs right the determinant is a {\sl graded line
bundle}, i.e., a pair $(\ell,\sL)$ where $\ell$ is an integer and $\sL$ a line
bundle. This is then used in the commutativity isomorphism $L\tensor M \riso
M\tensor L$ where the sign $(-1)^{\ell m}$ is used, where $\ell$ and $m$ are the
degrees of $L$ and $M$ respectively. The coherence conditions proved in
\cite{knudsen76} then show that we get an unambiguous isomorphism between
tensor products of the same graded line bundles in different order.

We shall need one property of the determinant beyond those of
\cite[Def.~4]{knudsen76}: If, for a perfect complex $C$ of $\sO_S$-modules, $S$
a scheme, we let $C^*:={\rm RHom}_S(C,\sO_S)$, then we have a canonical isomorphism
$\rho_C: \det(C)^* \iso \det(C^*)$ functorial for
quasi-isomorphisms. Indeed, this can be shown by verifying that $C \mapsto
(\det(C^*))^*$ verifies the conditions of \cite[Def.~4]{knudsen76} and hence by
\cite[Thm.~2]{knudsen76} it is (canonically) isomorphic to $C \mapsto \det
C$. (It can also be done by direct verification.) In any case note that we
define the dual of a graded line bundle as $(\ell,\sL)^*=(-\ell,\sL^*)$ and that we
identify $(\sL\tensor\sM)^*=\sM^*\tensor\sL^*$ by the pairing
\begin{displaymath}
(\sL\tensor\sM)\tensor(\sM^*\tensor\sL^*)=\sL\tensor(\sM\tensor\sM^*)\tensor\sL^*
\xrightarrow{{\rm id}\tensor {\rm ev}_{\sM}\tensor{\rm id}}
\sL\tensor\sO_S\tensor\sL^*=\sL\tensor\sL^* \to \sO_S,
\end{displaymath}
where we have used the above sign rule for the permutation. The unicity (as
well as direct computation) also gives that we have a commutative diagram

\begin{equation}
\label{determinant symmetry}
\begin{CD}
\det(C^*)^* @>\rho_C^*>> \det(C)^{**}\\
@VV\rho_{C^*}V           @AA {\rm ev}_{\det(C)}A\\
\det(C^{**})@<\det({\rm ev}_C)<< \det(C),
\end{CD}
\end{equation}
where ${\rm ev}_C: C \to C^{**}$ is the evaluation map (and similarly for
${\rm ev}_{\det(C)}$). If $\to \sE \to \sF \to \sG \to$ is a distinguished
triangle of perfect complexes we have a distinguished triangle $\to \sG^* \to
\sF^* \to \sE^* \to$ and the resulting identification
\begin{displaymath}
(\det\sE\tensor\det\sG)^*=\det(\sF)^*=\det(\sF^*)=\det(\sG^*)\tensor\det(\sE^*)=\det(\sG)^*\tensor\det(\sE)^*
\end{displaymath}
is then a special case of the above identification.

Now, let $\pi: X \to S$ be a smooth and proper map of schemes of pure relative
dimension $n$ over a base $S$ of positive characteristic $p \neq 2$. Let $\cL$ be
the determinant of $R\pi_*\Omega^\bullet_{X/S}$ (which exists as
$R\pi_*\Omega^\bullet_{X/S}$ is a perfect complex). By Poincar\'e duality we
have a canonical isomorphism $(R\pi_*\Omega^\bullet_{X/S})^* \riso
R\pi_*\Omega^\bullet_{X/S}[-2n]$ which upon taking determinants gives an
isomorphism
\begin{displaymath}
\cL^* \mapright{\rho} \det(R\pi_*\Omega^\bullet_{X/S}[-2n]) \riso \cL^{(-1)^{2n}}=\cL.
\end{displaymath}
Now, Poincar\'e duality gives a symmetric pairing and by (\ref{determinant
  symmetry}) this gives a perfect symmetric pairing $\cL\tensor\cL \riso \sO_S$.
On the other hand, the naive truncations\footnote{The reader who feels the need
  to recall the definition of naive and canonical truncations could profitably
  consult \cite{illusie04::topic}} of the de Rham complex give rise to
distinguished triangles
\begin{displaymath}
\to R\pi_*\Omega^{\geq i+1}_{X/S} \to R\pi_*\Omega^{\geq i}_{X/S} \to R\pi_*\Omega^i_{X/S}[-i]\to.
\end{displaymath}
Taking determinants we get (cf., \cite[Remark after Thm.~2]{knudsen76} for an
explication) an isomorphism

\begin{equation}\label{Hodge expansion}
\cL \riso \tensor_{i=0}^n\det(R\pi_*\Omega^i_{X/S})^{(-1)^i}.
\end{equation}
Similarly, we may use the canonical truncations to get a distinguished triangle
\begin{displaymath}
\to R\pi_*\sH^i(\Omega^\bullet_{X/S})[-i] \to R\pi_*\tau_{\geq i}\Omega^\bullet_{X/S}
\to R\pi_*\tau_{\geq i+1}\Omega^\bullet_{X/S}\to.
\end{displaymath}
Recall (cf., \cite[\S2.1]{illusie79::compl+rham}) the Cartier isomorphism
$\sH^i(F_{X/S*}\Omega^\bullet_{X/S})=\Omega^i_{X^{(p)}/S}$, where
$F_{X/S}: X \to X^{(p)}$ is the relative Frobenius map
fitting in the commutative diagram with Cartesian square and $F_X=W F_{X/S}$
\begin{displaymath}
\begin{xy}
\xymatrix{ X \ar[r]^{F_{X/S}} \ar[dr]_{\pi} &  X^{(p)} \ar[d]^{\pi^{(p)}} \ar[r]^W & X\ar[d]  \\
& S \ar[r]^{F_S} & S \\
}
\end{xy}
\end{displaymath}
where we will abuse notation and write $F_S$ for $W$.
Applying
$R\pi^{(p)}_*$, with $\pi^{(p)}: X^{(p)}\to S$ the structure map, we get
$R\pi_*\sH^i(\Omega^\bullet_{X/S})=R\pi^{(p)}_{*}\Omega^i_{X^{(p)}/S}
=R\pi^{(p)}_*F^*_S\Omega^i_{X/S}$.  Finally, using the base change formula
$R\pi^{(p)}_*LF^*_S=F^*_SR\pi_*$ we get an identification of
derived functors
$$
LF^*_SR\pi_*\Omega^i_{X/S}\riso R\pi_*\sH^i(\Omega^\bullet_{X/S}).
$$ 
Combining these formulas and taking determinants we obtain an isomorphism
(of graded line bundles)
\begin{equation}\label{conjugate expansion}
\cL \riso \tensor_{i=n}^0\det(LF^*_SR\pi_*\Omega^i_{X/S})^{(-1)^i} =
F^*_S\left(\tensor_{i=n}^0\det(R\pi_*\Omega^i_{X/S})^{(-1)^i}\right)
\end{equation}
where we note the inverse order due to the Cartier isomorphism.
We may then permute the last tensor product to get an isomorphism 
\begin{displaymath}
F^*_S\left(\Tensor_{i=n}^0\det(R\pi_*\Omega^i_{X/S})^{(-1)^i}\right)\riso 
F^*_S\left(\Tensor_{i=0}^n\det(R\pi_*\Omega^i_{X/S})^{(-1)^i}\right).
\end{displaymath}
Comparing the obtained formulas for $\cL$ and $F^*_S\cL$ we get an isomorphism
$\varphi: \cL \riso F^*_S\cL$; sometimes this is called an
{\sl $F$-structure} on $\cL$. Now, the isomorphism of (\ref{Hodge
expansion}) is compatible with duality (if we use the tensor product 
of the Serre
duality isomorphisms on the right hand side) and so is (\ref{conjugate
expansion}) because the Cartier isomorphism is multiplicative. This implies that
$\varphi$ is compatible with the pairing on $\cL$. We may now consider the sheaf
$L$ (in the \'etale topology on $S$) of fixed points under $\varphi$ and we know
that $\cL=L\tensor_{{\bf F}_p}\sO_S$, so that in particular $L$ is a local system of
$1$-dimensional ${\bf F}_p$-vector spaces. The pairing on $\cL$ induces a symmetric
non-degenerate pairing on $L$ and taking (locally) its discriminant gives us a
locally constant function from $S$ to ${\bf F}^*_p/{\bf F}^{*2}_p$. 
The latter group can
be identified, using the Legendre symbol, $\Legendre{-}{p}$, with $\{\pm1\}$ and
we shall call the resulting function $S \to \{\pm1\}$ the {\sl Hodge
discriminant} of $X \to S$. It clearly commutes with base change. In particular
its value can be computed fibrewise.

The Hodge discriminant uses the whole cohomology of $X/S$. Often we can also
work with the middle cohomology only.  Indeed, if we have an orthogonal or
symplectic $F$-zip $(M,C^\bullet,D_\bullet,\varphi_\bullet)$, we can make the same
construction, by making use of the two filtrations to identify the determinant
of the middle cohomology in two ways and compare them by $\varphi$: we
identify $\det M$ on the one hand with $\det\gr^\bullet C^\bullet$ and
on the other with 
$\det\gr_\bullet D_\bullet$, then use $\varphi$ to identify 
$F^*(\det\gr^\bullet C^\bullet)$ with $\det\gr_\bullet D_\bullet$ 
and finally use the induced
pairing to define a discriminant for the fixed points. This yields a Hodge
discriminant for the $F$-zip $(M,C^\bullet,D_\bullet,\varphi_\bullet)$.

Recall now (cf., \cite[Def.\ 3.1]{ogus82::hodge}) the definition of Ogus'
crystalline discriminant: Given an orthogonal or symplectic F-crystal $M$ over
an algebraically field ${\bf k}$ we get an induced F-crystal structure on $\det M$
for which $F$ is a power of $p$, $p^m$ say, times an isomorphism. Dividing by
$p^m$ and taking fixed points we get a ${\bf Z}_p$-module of rank $1$ with a perfect
pairing. Taking its discriminant and reducing modulo $p$ gives us an element of
${\bf F}_p^*/{\bf F}_p^{*2}$, the {\sl crystalline discriminant}. 

Before stating how these discriminants are related we recall that for a
proper smooth variety of pure dimension $n$ over a field ${\bf k}$ of 
positive characteristic $p$ we have the $\ell$-adic Betti number
$b_n(X)$ (which is the same as the rank of the $n$'th 
crystalline cohomology group), 
the number $b_n^{\prime}(X):=\dim H^n_{dR}(X/{\bf k})$ 
and the Hodge numbers $h^{ij}$ satisfying 
$$
b_n(X)\leq b_n^{\prime}(X) \leq \sum_{i+j=n} h^{ij}(X).
$$
If $b_n^{\prime}(X) =\sum_{i+j=n} h^{ij}(X)$ then the 
$E_2^{i,j}$-term of the Hodge-to-de Rham spectral sequence equals 
the $E_{\infty}^{i,j}$-term for all $i+j=n$. 
By dimension counting this then implies
the same thing for the $E_1^{i,j}$-term of the conjugate spectral
sequence. Hence the Hodge and conjugate filtrations on $H^n_{dR}(X/{\bf k})$ together
with the Cartier isomorphisms give an $F$-zip structure on $H^n_{dR}(X/{\bf k})$.
Cup product induces a non-degenerate pairing.
This is symplectic if $n$ is odd and orthogonal for $n$ even.
\begin{proposition}\label{Hodge discriminant}
Suppose that $X$ is a smooth and proper variety of pure dimension $n$ over a field
${\bf k}$ of positive characteristic $p$.

\item{i)} Assume that $b'_n:=\dim_{\bf k} H^{n}_{dR}(X/{\bf k})=\sum_{i+j=n}h^{ij}(X)$.
The Hodge discriminant of $X$ is equal to $\Legendre{-1}{p}^{(c_2-b'_n)/2}$
times the Hodge discriminant of the $F$-zip $H^n_{dR}(X/{\bf k})$, where $c_2$ is the
crystalline (=\'etale) Euler characteristic of $X$. If $n$ is odd the Hodge
discriminant is equal to $(-1)^{c_2/2}$.

\item{ii)} If $b_n(X)=\sum_{i+j=n}h^{ij}(X)$, then the Hodge discriminant of the
$F$-zip $H^n_{dR}(X/{\bf k})$ is equal to Legendre symbol of the crystalline
discriminant of the F-crystal $H^n_{cris}(X/{\bf W})$.
\begin{proof}
We start with the easily proven fact that if $\to X^\cdot \to
Y^\cdot \to Z^\cdot \to X^\cdot[1] \to$ is a distinguished triangle of complexes
(over a field) such the induced map $H^n(Z^\cdot) \to H^{n+1}(X^\cdot)$ is zero,
then we get a diagram, all of whose rows and columns are distinguished,
\begin{displaymath}
\begin{CD}
@>>>  \tau_{\leq n}X^\cdot @>>> \tau_{\leq n}Y^\cdot @>>>  \tau_{\leq n}Z^\cdot @>>> (\tau_{\leq n}X^\cdot)[1]@>>>\\
@.@VVV @VVV @VVV @VVV \\
@>>>  X^\cdot @>>> Y^\cdot @>>> Z^\cdot @>>> X^\cdot[1] @>>>\\
@.@VVV @VVV @VVV @VVV \\
@>>>  \tau_{> n}X^\cdot @>>> \tau_{> n}Y^\cdot @>>>  \tau_{> n}Z^\cdot @>>> (\tau_{> n}X^\cdot)[1]@>>>
\end{CD}
\end{displaymath}
Note furthermore that it is one of the properties of the Knudsen-Mumford
determinant (cf., \cite[Def~4]{knudsen76}) that the two ways of using this
diagram to get an isomorphism
\begin{equation}\label{3x3 compatability}
\det(Y^\cdot)\riso \det(\tau_{\leq n}X^\cdot)\Tensor\det(\tau_{>
  n}X^\cdot)\Tensor \det(\tau_{\leq n}Z^\cdot)\Tensor \det(\tau_{> n}Z^\cdot)
\end{equation}
give the same result.

The degeneration of the Hodge to de Rham spectral sequence and that of the
conjugate spectral sequence both at total degree $i+j=n$ implies that the
necessary conditions are fulfilled to apply this and thus allow us to
conclude that we have distinguished triangles:
\begin{eqnarray*}
\to \tau_{<n}R\Gamma(X,\Omega^{\geq i+1}) \to  \tau_{<n}R\Gamma(X,\Omega^{\geq i})\to
\tau_{<n}R\Gamma(X,\Omega^{i}[-i])\to\\
\to \tau_{>n}R\Gamma(X,\Omega^{\geq i+1}) \to  \tau_{>n}R\Gamma(X,\Omega^{\geq i})\to
\tau_{>n}R\Gamma(X,\Omega^{i}[-i])\to\\
\to \tau_{<n}R\Gamma(X,\sH^i[-i]) \to  \tau_{<n}R\Gamma(X,\tau_{\geq i}\Omega^{\bullet})\to
\tau_{<n}R\Gamma(X,\tau_{\geq i+1}\Omega^{\bullet})\to\\
\to \tau_{>n}R\Gamma(X,\sH^i[-i]) \to  \tau_{>n}R\Gamma(X,\tau_{\geq i}\Omega^{\bullet})\to
\tau_{>n}R\Gamma(X,\tau_{\geq i+1}\Omega^{\bullet})\to\\
\end{eqnarray*}

This gives us expansions
\begin{eqnarray*}
\det(\tau_{<n}R\Gamma(X,\Omega^\bullet))&=&\Tensor_{i=0}^n\det(\tau_{< n-i}R\Gamma(X,\Omega^i))^{(-1)^i}\\
\det(\tau_{>n}R\Gamma(X,\Omega^\bullet))&=&\Tensor_{i=0}^n\det(\tau_{> n-i}R\Gamma(X,\Omega^i))^{(-1)^i}\\
\det(H^n(X,\Omega^\bullet))&=&\Tensor_{i=0}^n\det(H^{n-i}(X,\Omega^i))^{(-1)^i}
\end{eqnarray*}
and
\begin{eqnarray*}
\det(\tau_{<n}R\Gamma(X,\Omega^\bullet))&=&\Tensor_{i=n}^0\det(\tau_{< n-i}F^*R\Gamma(X,\Omega^i))^{(-1)^i}\\
\det(\tau_{>n}R\Gamma(X,\Omega^\bullet))&=&\Tensor_{i=n}^0\det(\tau_{> n-i}F^*R\Gamma(X,\Omega^i))^{(-1)^i}\\
\det(H^n(X,\Omega^\bullet))&=&\Tensor_{i=n}^0\det(F^*H^{n-i}(X,\Omega^i))^{(-1)^i},
\end{eqnarray*}
where $F=F_{{\rm Spec}({\bf k})}$.

Now, we also have an expansion
\begin{equation}\label{n-decomposition}
\det(R\Gamma(X,\Omega^\bullet))=
\det(\tau_{<n}R\Gamma(X,\Omega^\bullet))\Tensor\det(H^n(X,\Omega^\bullet))\Tensor
\det(\tau_{>n}R\Gamma(X,\Omega^\bullet)).
\end{equation}
We have already given the left hand side an $F$-structure and the isomorphisms
above give $F$-structures on each factor of the right hand side. However,
it follows from (\ref{3x3 compatability}) and the fact that the tensor product
on graded line bundles is symmetric monoidal (see \cite[Ch.\ I]{knudsen76})
that those two $F$-structures are the same.

As for the self-pairing, the duality induces isomorphisms
$(\tau_{<n}R\Gamma(X,\Omega^\bullet_X))^*=\tau_{>n}R\Gamma(X,\Omega^\bullet_X)[-2n]$
and $H^n(X,\Omega^\bullet_X)^*=H^n(X,\Omega^\bullet_X)$. This means that in the
decomposition of (\ref{n-decomposition}) the self-pairing on the left
corresponds to a pairing on the right which pairs the first factor to the third
and the second to itself. This is compatible with semi-linear structure so that
when we take fixed points under $F$ we get a decomposition $L=L^<\Tensor
L^=\Tensor L^>$, with $L^<$ and $L^>$ being paired to each other and $L^=$ to
itself. Taking the signs into account we get that the discriminant of $L$ is
equal to $\Legendre{-1}{p}^{\dim L^<}$ times the discriminant of $L^=$. However,
the dimension of $L^<$ is equal to the dimension of
$\tau_{>n}R\Gamma(X,\Omega^\bullet_X)$ which is $(c_2-b'_n)/2$. Of course, when
$n$ is odd $L^=$ is trivial. This concludes the proof of i).

As for ii) we are reduced by i) to showing that the Hodge
discriminant of $H^n_{dR}(X/{\bf k})$ is equal to the crystalline discriminant of
$H^n_{cris}(X/{\bf W})$. By the Mazur-Ogus result \cite[Appendix]{berthelot78::notes}
we may choose a basis $e_1^{1},\dots,e_{h^1}^1,e_1^2,\dots,e_{h^r}^r$ of
$H^n_{cris}(X/{\bf W})$ such that $Fe_j^s$ is divisible by $p^s$ and such that the
Hodge filtration of $H^n_{dR}(X/{\bf k})$ is given by $H_i=\sum_{j\leq i}\sum_{1\leq
  k \leq h^j}{\bf k} \overline{e}_k^{j}$, the conjugate filtration is given by
$H^c_i=\sum_{j\leq i}\sum_{1\leq k \leq h^j}{\bf k} \overline{p^{-j}Fe}_k^{j}$ and
the inverse Cartier isomorphism is given by $\{p^{-j}F\}$. Unraveling
definitions then gives part ii).
\end{proof}
\end{proposition}

\begin{remark}
  As the proposition shows, under suitable conditions the Hodge discriminant is
  equivalent to the crystalline discriminant. We justify introducing new
  notation since it would be somewhat artificial to use ``crystalline'' in a
  situation where it is not relevant; moreover, it makes sense more generally,
  for instance in the case of Enriques surfaces in characteristic two, cf.\
\cite{EHS}. Please
  note that we have defined the Hodge discriminant to be the Legendre symbol
  applied to an element of ${\bf F}^*_p/{\bf F}^{*2}_p$ rather than the element itself.
  (Of course the Legendre symbol gives a bijection with this group and
  $\{\pm1\}$ so no information is lost.) The reason for this convention is to
  make the formulas of Proposition \ref{Odd Hodge discriminant} as nice as
  possible; otherwise that formula would have to involve the inverse of the
  isomorphism provided by the Legendre symbol.
\end{remark}
The further properties of the Hodge discriminant will differ somewhat depending on
whether $H$ has even or odd dimension, so we shall discuss each case separately.
\begin{subsection}{The Hodge discriminant of an orthogonal flagged $F$-zip}

We now derive a formula for the Hodge discriminant of an orthogonal
flagged $F$-zip. In the odd-dimensional case we need one more notion.  Hence
consider a flagged orthogonal $F$-zip $(H,C^\bullet,D_\bullet,\varphi_\bullet)$
of dimension $2m+1$. We get two induced isomorphisms
\begin{displaymath}
C^{m-1}/C^{m} \liso C^{m-1}\cap D_{m+1}/C^m\cap D_m \riso D_{m+1}/D_{m}
\end{displaymath}
and together with the inverse Cartier isomorphism they give rise to an
isomorphism
\begin{displaymath}
F^*(C^{m-1}/C^{m}) \riso C^{m-1}/C^{m}.
\end{displaymath}
On the other hand, we also have a pairing on $C^{m-1}/C^{m}$ induced from that
of $H$ and it is compatible with $F^*(C^{m-1}/C^{m}) \riso
C^{m-1}/C^{m}$. Hence, we get an ${\bf F}_p^*/{\bf F}_p^{2*}$-valued discriminant by
taking fixed points. We shall call it the {\sl middle discriminant}.
\begin{proposition}\label{Odd Hodge discriminant}
Let $(H,C^\bullet,D_\bullet,\varphi_\bullet)$ be an orthogonal flagged $F$-zip.

\item{i)} Assume $H$ has dimension $2m+1$ of type $w \in W^C_m$. Then
the Hodge discriminant equals $(-1)^{n_s}\Legendre{d}{p} {\disc}{(w)}$, where
$d=(-1)^md'$ with $d'$ the middle discriminant, and $s=[r/2]$.
\item{ii)} Assume $H$ has dimension $2m$ of type $w \in W'^D_m$. Then the Hodge
discriminant equals $(-1)^{n_s}{\Legendre{-1}{p}}^m{\disc}{(w)}$, where $s=[r/2]$.
\begin{proof}
  For the odd case we may assume by Proposition \ref{Versal $F$-zips} that the
  $F$-zip is associated to a ${\bf k}$-point of $Y_w$, i.e., there is a basis
  $e_1,\dots,e_{2m+1}$ for $H$ with $\langle e_i, e_j \rangle = \delta_{i,2m+2-j}$,
  $F_i=\sum_{j\leq i}{\bf k} e_{j}$, and $D_i=\sum_{j\leq i}{\bf k} e_{w^{-1}(j)}$ with
  $\varphi_k$ acting as specified by the construction of the $F$-zip on $Y_w$.
  This implies that $\varphi_k(e_{n_k+1}\land e_{n_k+2}\land\cdots\land
  e_{n_{k+1}})=\epsilon_ke_{w^{-1}w_\emptyset(n_k+1)}\land
  e_{w^{-1}w_\emptyset(n_k+2)}\land\cdots\land e_{w^{-1}w_\emptyset(n_{k+1})}$.
  Here $\epsilon_k=1$ when $k\ne m$ and $\epsilon_m=\pm1$ with $+1$ when the
  middle discriminant is a square and $-1$ when it is not. This implies that the
  semi-linear map on the determinant takes $e_1\land e_2\land\cdots\land
  e_{2m+1}$ to $\epsilon_me_{w^{-1}w_\emptyset(1)}\land
  e_{w^{-1}w_\emptyset(2)}\land\cdots\land
  e_{w^{-1}w_\emptyset(2m+1)}=\epsilon_m{\disc}(w^{-1}w_\emptyset) e_1\land
  e_2\land\cdots\land e_{2m+1}$. Similarly we have that $\langle e_1\land
  e_2\land\cdots\land e_{2m+1}, e_1\land e_2\land\cdots\land
  e_{2m+1}\rangle =(-1)^m$. Now we conclude by
the mod $p$ version of \cite[Formula 3.4]{ogus82::hodge}
using that $d'$ is a square precisely when
  is $\epsilon_m$ and that ${\disc}(w_\emptyset)=(-1)^{n_s}$.

  The proof of the even-dimensional case is identical to the odd case except
  that we always have $\epsilon_m=1$.
\end{proof}
\end{proposition}

\end{subsection}
\end{section}

\begin{section}{K3 Surfaces}

In this section we shall consider the primitive cohomology of a polarized 
K3 surface and show that it possesses a minimal {\sl stable} filtration
that refines the Hodge filtration and that it can be refined to a 
so-called {\sl final} filtration if the field of definition is separably closed.

The results of the preceding section can be applied because the crystalline
cohomology of a K3 surface is without torsion and the Hodge-to-de Rham
spectral sequence degenerates  at the $E_2$-level, cf.\ \cite{De81,illusie79::compl+rham}. 

Recall that for a K3 surface the N\'eron-Severi
group $\NS(X)$ is equal to the Picard group of $X$.
Let $N$ be a non-degenerate integral lattice. A {\sl (partial)
$N$-marking} of a K3 surface $X$ over a field ${\bf k}$ of positive
characteristic $p$ is an isometric embedding $N \to \NS(X)$.
The {\sl discriminant} of the marking
is the discriminant of the lattice $N$.  We shall only be interested
in partial markings whose degree (order of the discriminant group) 
is prime to $p$ and thus that will be assumed
unless otherwise mentioned. We define the {\sl primitive cohomology} of
an $N$-polarized K3 surface $X$ as the orthogonal complement of the 
image of $N$ in $H^2_{dR}(X/{\bf k})$.

The primitive cohomology is an orthogonal $F$-zip with dimension vector for its
Hodge filtration being $(0,1,n-1,n)$ for some $n$. We shall also need another
type of $F$-zip. Namely, an $F$-zip of dimension 
type $(0,\dots,0,m,\dots,m)$ shall
be called a {\sl Tate $F$-zip}. Tate $F$-zips thus consist of an
orthogonal vector space $V$ and an orthogonal $F$-structure
$\varphi: {F^*V} \to {V}$. It is thus completely described by the orthogonal
representation of the Galois group of ${\bf k}$ given by the action on
$\sV:=\ker(\varphi-1)$ on $V\Tensor_{\bf k}\overline{\bf k}$. 
We shall say that the Tate $F$-zip is {\sl split} resp.\ {\sl non-split} 
as the form on $\sV$ is. Its Hodge discriminant is then  
$\Legendre{d}{p}$, where $d$ is
the discriminant of $\sV$. In these terms we have that $H^2_{dR}(X/{\bf k})$ 
is the sum, as $F$-zip, of the primitive cohomology and $N\Tensor{\bf k}$ 
considered as a Tate $F$-zip.
\begin{definition}
Let $M$ be a stable partially flagged orthogonal or symplectic $F$-zip.
\begin{itemize}
\item[i)] $M$ is {\sl final} if it is complete.

\item[ii)] If $M$ is symplectic or orthogonal of odd dimension then it is
{\sl canonical} if every stable flag is a refinement of it.

\item[iii)] If $M$ is orthogonal of even dimension it is {\sl canonical} if
every stable flag is a refinement of $M$ or possibly, when it exists, its flip.
\end{itemize}
\end{definition}
\begin{example}
For a Tate $F$-zip its (trivial) Hodge filtration already is canonical. A final
filtration is an ${\FF}_p$-rational self-dual flag except in case the
${\FF}_p$-form is non-split. Then the middle element of the flag is only defined over
${\FF}_{p^2}$.
\end{example}
\begin{lemma}\label{Final r-characterisation}
Let $w$ be an element of $\Weyl{B}{m}$ or $\Weylp{D}{m}$. Assume that for all $1\leq
i,j \leq n-1$, where we do not have $i=j=n/2$,
\begin{displaymath}
r_w(i,j) =
\begin{cases}
\min(j,r_w(i,n-1)+1)-1 & \text{ if }i<a,\\
\min(j,r_w(i,n-1)) & \text{ if }i\geq a,
\end{cases}
\end{displaymath}
where $a:=w^{-1}(1)$ and $n=2m+1$ in the $B$-case and $2m$ in the $D$-case. Then $w$
is a final element and conversely the $r_w$ for $w$ final fulfils this condition.
\begin{proof}
By definition we have 
$r_w(i,n-1)=\#\{ 1\leq b\leq i: w(b)\leq n-1 \}$. Thus it is clear that it 
is determined completely by $w^{-1}(n)$ which is equal to $n+1-a$. 
The assumed conditions on
$r_w$ then implies that the whole function is determined by $a$ and hence so is
$w$. It is easy to verify that a final $w$ fulfils the conditions and that there
is one such element for each $a$.
\end{proof}
\end{lemma}
Recall that a flagged $F$-zip has a type which is an element $w$ of a Weyl group.
We now link the definition of a final $F$-zip to the notion of a final element
of a Weyl group.

\begin{proposition}
A flagged orthogonal $F$-zip of type $(0,1,n-1,n)$ is final precisely when its
type is a (twisted) final element.
\begin{proof}
Suppose that $E_{\bullet}$ is a final filtration and $G_{\bullet}$ the
corresponding conjugate filtration (so that $0\subset E_1 \subset E_{n-1}
\subset E_n$ is the Hodge filtration and $0\subset G_1 \subset G_{n-1}
\subset G_n$ the conjugate filtration). For each $1 \leq i\leq n-1$ we have by
assumption that for every $i$ the subspace $G_i\cap E_{n-1}+E_1$ is equal
to some $E_r$ (or possibly its flip $E'_r$ if $2r=n$). Then for $1\leq j\leq
n-1$
\begin{displaymath}
G_i\cap E_j+E_1 = (G_i\cap E_{n-1}+E_1)\cap E_j = E_r\cap E_j = E_{\min(r,j)},
\end{displaymath}
where the end result would instead be $E_{r-1}$ if $G_i\cap
E_{n-1}+E_1=E'_r$ and $j=r$. Now, $r_w(i,j)=\dim(G_i\cap E_j)$ and in
particular $E_1\subseteq G_i$ precisely when $i\geq w^{-1}(1)$ and thus
$\dim(G_i\cap E_j+E_1)$ is equal to $r_{w}(i,j)+1$ if $i<w^{-1}(1)$ and
$r_{w}(i,j)$ otherwise (supposing that we do not have $i=j=n/2$). This shows
that $r_w$ fulfils the condition of Lemma \ref{Final r-characterisation} and
hence $w$ is final. The converse is just a matter of tracing the argument
backwards.
\end{proof}
\end{proposition}
Recall that two orthogonal $F$-zips of type $(0,1,n-1,n)$
are called {\sl opposite} if their intersections
have the smallest possible dimensions, i.e., $F_1 \not\subset E_{n-1}$.
It follows from either the description of the final elements or from the proof
of the next theorem that the canonical filtrations have the form $U_1\subset
U_2\subset\cdots\subset U_k\subset U_{n-k}\subset\cdots \subset U_n$, where the
primitive cohomology has dimension $n$. We shall call $U_{n-k}/U_k$ the
{\sl middle part} of the canonical filtration. It comes equipped with a
quadratic form induced from that of $H^2_{dR}(X/{\bf k})$ 
and the Cartier isomorphism
induces an orthogonal $p$-linear isomorphism of it (i.e., a Tate $F$-zip
structure). The fixed points under the Cartier isomorphism then give an
${\FF}_p$-rational structure on the middle part and the quadratic form induces a
quadratic form on it. We shall say that the canonical filtration is
{\sl split} resp.\ {\sl non-split} according to as that form is.
\begin{theorem}\label{Final/canonical filtrations}
Let $X$ be a polarised K3 surface of degree prime to $p$ over a field ${\bf k}$
 of characteristic $p >0$ and let $H$ be its primitive Hodge cohomology
of dimension $n$ with $m:=[n/2]$. Then $H$ has a canonical filtration.
Any final filtration is obtained from the canonical one by choosing a
complete $F$-stable filtration. If ${\bf k}$ is separably closed
$H$ has a final filtration. All final filtrations have the same (twisted) final
type.
\begin{proof}
  We start with the induced Hodge filtration $0 \subset E_1 \subset E_{n-1}
  \subset E_n=H$ on the primitive cohomology with conjugate filtration $0
  \subset F_1 \subset F_{n-1}\subset F_n=H$ with $F_i=E_{n-i}^c$.  If $F_1=E_1$
  then the two filtrations coincide and then this partially flagged $F$-zip is
  canonical as one easily checks. If $F_1 \not= E_1$ then we consider the image
  of $F_1$ in $E_n/E_{n-1}$. If this image is non-zero then the Hodge filtration
  and the conjugate one are opposite and we get a stable and hence canonical
  flagged $F$-zip.  So suppose that $F_1$ has non-zero image in $E_{n-1}/E_1$.
  We can apply Frobenius and use the Cartier isomorphism to get the image
  $\overline{F}_2$ in $E_{n-1}^c/E_1^c=F_{n-1}/F_1$. We then add to our flag the  inverse image $F_2$ of $\overline{F}_2$ in $F_{n-1}$.  Now $F_1$ is totally
  isotropic, hence its image in $E_{n-1}/E_1$ is as well and as the Cartier
  isomorphism is multiplicative for the wedge product, so is $F_2$ and therefore
  $F_{n-2}:=F_2^{\perp}$ contains $F_2$.  We then continue this process: if
  $F_2$ is not contained in $E_{n-1}$ then we have obtained a stable filtration.
  This is also the case if $E_1\subset F_2$. On the other hand if $F_2 \subset
  E_{n-1}$ and $F_2\cap E_1 = \{0\}$ we consider the image of $F_2$ in
  $E_{n-1}/E_1$ and transfer via $F$ and the Cartier isomorphism the image to
  $F_{n-1}/F_1$. Note that each stage of the induction $F_{n-i} \subset E_{n-1}$
  precisely when $E_1\subset F_i$ and thus we will not be forced to introduce
  any new elements to the flag because of the position of $F_{n-i}:=F_i^\perp$.
  It is clear that this process stops and gives a canonical flag.  That a
  canonical filtration can be extended to a final filtration and that they are
  all of the same type is easy and similar to the abelian case, see
  \cite[Def-Lemma 2.11]{ekedahl10::cycle+class+e+o+strat}.
\end{proof}
\end{theorem}

\end{section}

\begin{section}{Canonical filtrations versus the height and Artin invariant}
As we just saw, the primitive part of the 2nd de Rham cohomology of a 
K3 surface in positive characteristic comes with a canonical filtration. 
If our field is separably closed we can refine it to a final filtration.
  A natural question is now what the type of the final filtration means
  geometrically. The following theorem will provide the answer.
It relates the relative position of a final filtration and its conjugate one, 
given by an element $w$ of a Weyl group, to geometric invariants.
Recall that we
  have two invariants for a K3 surface $X$ in positive characteristic, the
  {\sl height} and if the height is infinite we also have 
the so-called {\sl Artin invariant}.  The height $h(X)$ is
  the height of the formal Brauer group, a smooth formal group of dimension $1$.
  This invariant assumes values $1\leq h \leq 10$ or $h=\infty$; in the latter
  case the formal Brauer group is the formal additive group. The Artin invariant
  $\sigma_0$ can be defined for supersingular K3 surfaces, i.e.\ 
those with   $h=\infty$ by putting 
${\rm disc}(H^2(X,\ZZ_p(1)))=-p^{2\sigma_0}$, cf.,
  \cite{artin74::super+k3}.  We then have $1 \leq \sigma_0 \leq 10$.  The case
  $h=1$ is called the ordinary case and it is 
the generic finite height case and $\sigma_0=10$ is the generic
  supersingular case.

These invariants can be detected by the
crystalline cohomology.  We therefore recall first some facts on crystalline
cohomology and the relation with de Rham-cohomology (cf., \cite[Thm
8.26]{berthelot78::notes}).

Let ${\bf W}({\bf k})$ be the ring of Witt vectors of ${\bf k}$ and 
let $\sigma$ be the map
on ${\bf W}({\bf k})$ induced by the Frobenius map on ${\bf k}$.  
The second crystalline cohomology group ${\mathcal H}:=H^2(X/{\bf W}({\bf k}))$ is a free ${\bf W}({\bf k})$-module of rank
$22$ and is provided with a ${\bf W}({\bf k})$-linear map 
$F: \sigma^*{\cH}\to \cH$.
We have a natural isomorphism from ${\cH}/p{\cH} \riso H^2_{dR}(X/{\bf k})$ 
and by base change by the Frobenius map on ${\bf k}$ we get an isomorphism
$\sigma^*{\cH}/p\sigma^*{\cH} \riso H^2_{dR}(X^{(p)}/{\bf k})$.  
If we put ${\cH}_i:= F^{-1}p^{2-i}{\cH}$ for $i=0,1,2$, 
then the images $H_i$ of the ${\cH}_i$ in
$\sigma^*{\cH}/p\sigma^*{\cH}$ give the Hodge filtration on
$H^2_{dR}(X^{(p)}/{\bf k})$ 
(with $E_1=H_0$, $E_{n-1}=H_1$ and $E_n=H_2$ in the
notation of the proof of Thm \ref{Final/canonical filtrations}) while the images
of the ${\cH}^c_{i}:=p^{-i}F{\cH}_{2-i}$ in $H^2_{dR}(X/{\bf k})$ 
give the conjugate filtration. 
Finally, the inverse Cartier isomorphism is induced by the map
$p^{-i}F: {\cH}_{2-i}\to {\cH}^c_{i}$.

\smallskip
We now give the main result connecting the final type with the classical 
invariants (height and Artin invariant) and the Hodge discriminant.

Recall that we choose a marking by giving an isometric
embedding $N \to {\rm NS}(X)$. In particular we have a discriminant
$d$ of the marking.

\begin{theorem}\label{Final/canonical filtrations type}
Let $X$ be a polarized K3 surface of degree prime to $p$ over a field ${\bf k}$ of
characteristic $p >0$ and let $H$ be its primitive Hodge cohomology
of dimension $n$ with $m:=[n/2]$.
\begin{enumerate}
\item[i)] If $X$ has finite height $h$ with $2h < n$, then $H$ has final type
$w_{h}$ or $w'_{h}$. When $n$ is even it is $w_h$ if the middle part is
non-split and $w'_h$ if it is split.

\item[ii)] If $X$ has finite height $h = n/2$, then $H$ has final
type $w'_{m}$.

\item[iii] If $X$ is supersingular with Artin invariant $\sigma_0<n/2$, then $H$ has
final type $w_{2m+1-\sigma_0}$ or $w'_{2m+1-\sigma_0}$. When $n$ is even it is
$w_{2m+1-\sigma_0}$ if the middle part is split and $w'_{2m+1-\sigma_0}$ if it is
non-split.

\item[iv)] If $X$ is supersingular with Artin invariant $\sigma_0=n/2$, then
$H$ has final type $w_{m+1}$.

\item[v)] The Hodge discriminant of $H$ is equal to $\Legendre{-d}{p}$, where
$d$ is the discriminant of the marking.
\end{enumerate}
\begin{proof}

  Note that because the discriminant of the marking is prime to $p$, our space
  ${\cH}:=H^2(X/{\bf W}({\bf k}))$ splits into the orthogonal direct sum 
$({\bf W}({\bf k}) \Tensor N^{\perp})\oplus  ({\bf W}({\bf k})\Tensor N)$, 
where $N$ embeds using the crystalline Chern class and a
  similar statement is true for Hodge cohomology. This gives in particular that
  the Hodge discriminant of $H^2_{dR}(X/{\bf k})$ is the product of the Hodge
  discriminant of the primitive part and the Legendre symbol of the discriminant
  of $N/pN$. By a theorem of Bloch and Ogus (cf., \cite[Thm
  4.9]{ogus83::ae+torel+k3}), it is equal to $(-1)^{22-1}$ and this together
  with the relation between the Hodge and crystalline discriminants gives~v).

  Now, if we perform our construction of the canonical filtration on all of
  $H^2_{dR}(X/{\bf k})$, 
it will be performed separately on the reduction modulo $p$
  of ${\bf W}({\bf k})\Tensor N^\perp$ and ${\bf W}({\bf k})\Tensor N$. 
Furthermore it will be completely trivial
  on the second factor having a canonical flag consisting only of the zero
  subspace and the full space. Hence we may as well work with the full
  crystalline and de Rham cohomologies rather than their primitive parts and we
  shall do exactly that. Thus now $H$ may be identified with ${\cH}/p{\cH}$.
  With these results in mind we shall now consider the different cases.

{\sl Case 1.}
Consider first the case of finite height $h$. Then ${\cH}$ splits as an
orthogonal $F$-stable direct sum
$$
M_{1/h}\oplus {\bf W}({\bf k})^{22-2h}(1)\oplus M_{2-1/h}.
$$
Here $M_{1/h}$ is the crystalline Dieudonn\'e module with basis
$e_1,e_2,\dots,e_{h-1}$ where $Fe_i=p\, e_{i-1}$ for $i=2,\dots,h-1$ and
$Fe_1=e_{h-1}$. Further, ${\bf W}({\bf k})(1)^{22-h}$ 
is free of rank $22-2h$ with $F$
acting as $p$ on a basis. Finally, $M_{2-1/h}$ is the dual $M_{1/h}^*(1)$ of
$M_{1/h}$ as Dieudonn\'e module twisted once (i.e., the Frobenius map is
multiplied by $p$). In particular, $M_{2-1/h}$ has a basis
$f_1,f_2,\dots,f_{h-1}$ with $Ff_i=pf_{i+1}$ for $i=1,\dots,h-2$ and
$Ff_{h-1}=p^2f_1$. Furthermore, we have an orthogonal decomposition
$M_{1/h}\oplus M_{2-1/h}\perp {\bf W}({\bf k})^{22-2h}(1)$ 
which again means that the Hodge
and conjugate filtrations will be a direct sum of those of the summands. As
before the filtrations on the ${\bf W}({\bf k})^{22-2h}(1)$ 
factor will be trivial and we
hence may restrict to the other factor and will put ${\cH}$ equal to
$M_{1/h}\oplus M_{2-1/h}$.
There the pairing will be given by identifying $M_{2-1/h}$ with the
dual of $M_{1/h}$.
From the description above we conclude that (employing 
the notation $(M_{1/h})_i=M_{1/h} \cap{\mathcal H}_i$ for $i=0,1$ 
with ${\mathcal H}_i$ as defined above)
$$
\begin{aligned}
(M_{1/h})_1 &=p{\bf W} e_1+{\bf W} e_2+\cdots+{\bf W} e_{h-1},\\
(M_{2-1/h})_1 &=M_{2-1/h},\\ 
(M_{1/h})_0 &=p^2{\bf W} e_1+p{\bf W} e_2+\cdots+p{\bf W} e_{h-1},\\
(M_{2-1/h})_0 &=p{\bf W} f_1+\cdots+p{\bf W} f_{h-2}+{\bf W} f_{h-1}. \\
\end{aligned}
$$
This implies that $H_0$ has $\overline{f}_{h-1}$
as a basis and $H_1$ has
$
\overline{e}_2,\dots,\overline{e}_{h-1},\overline{f}_1,\dots,\overline{f}_{h-1}
$
as a basis. Similarly, we get that
$H^c_0$ has $\overline{e}_{h-1}$
as a basis and $H^c_1$ has
$\overline{e}_1,\dots,\overline{e}_{h-1},\overline{f}_1,\dots,\overline{f}_{h-2}$
as a basis and we also see that
$C^{-1}: F^*(H_1/H_0)\to H^c_1/H^c_0$
takes $\overline{e}_i$ to $\overline{e}_{i-1}$ for $1\leq i\leq h-1$ and
$\overline{f}_i$ to $\overline{f}_{i+1}$ for $0\leq i\leq h-2$.

Now, as we saw during the construction of the canonical filtration we do not
need to introduce $U_{n-i}:=U_i^\perp$ of our desired filtration at each stage
of the construction but can do it when the construction is finished. From the
description above it follows that
$H_0+H^c_0={\bf k}\overline{e}_{h-1}+{\bf k}\overline{f}_{h-1}$.  Transferring by the
Cartier isomorphism forces us to add ${\bf k}\overline{e}_{h-2}+{\bf k}
\overline{e}_{h-1}$
to the refinement of the conjugate filtration. Continuing one sees that all the
${\bf k}\overline{e}_{h-i}+\cdots+{\bf k}\overline{e}_{h-1}$ 
for $i\leq h$ must be added
to the canonical filtration. Then also their annihilators
${\bf k}\overline{f}_{n-i}+\cdots+{\bf k}\overline{f}_{h-1} +
\overline{M_{1/h}}$ must be
added. We thus get a canonical filtration with the property that
${\bf k}^{22-2h}=\overline{ {\bf W}({\bf k})^{22-2h}(1)}$ 
maps isomorphically to the quotient
of the $h$'th and $h+1$'st step in the filtration.  We can complete such a
canonical flag by adding a complete self-dual flag of
$\overline{{\bf W}({\bf k})^{22-2h}(1)}$ which is fixed under $p^{-1}F$ i.e., an
${\FF}_p$-rational such flag.  By comparing our zip to the standard case of
Proposition \ref{Flags classification} and the form of the final elements one
sees directly that these are of type $w_h$ or $w_h^{\prime}$. To decide whether
the form on ${\FF}_p^{22-2h}$ is split or not we interpret its discriminant in
terms of the crystalline discriminant (cf., \cite{ogus82::hodge}), i.e., the
discriminant of the fixed points of $p^{2h-22}F$ on
$\Lambda^{22-2h}({\bf W}({\bf k})^{22-2h}(1))$ 
multiplied by $(-1)^{11-h}$. As ${\cH}$
splits up as the orthogonal direct sum of ${\bf W}({\bf k})^{22-2h}(1)$ and a hyperbolic
space on $M_{1/h}$ we see that the crystalline discriminant of ${\cH}$
equals the product of $(-1)^h$ and the crystalline discriminant of
${\bf W}({\bf k})^{22-2h}(1)$. It follows from Proposition \ref{Odd Hodge discriminant}
that if the type is $w$, then the Hodge discriminant is
$-{\Legendre{-1}{p}}^{11}{\disc}(w)$. Now, from the Bloch-Ogus theorem and what we
just proved we get that the Hodge discriminant of the middle part is
$-{\Legendre{-1}{p}}^{h}{\disc}(w)$ and as it is split precisely when its Hodge
discriminant is ${\Legendre{-1}{p}}^{11-h}$ we see that it is split precisely
when ${\disc}(w)=-1$.

{\sl Case 2.}  Turning now to the case of infinite height let us recall the
setup of \cite{ogus79::super+k3}. (As we do not want to assume that $\rho=22$
we shall however replace $\NS\Tensor\ZZ_p$ by the flat cohomology group
$H^2(X,\ZZ_p(1))$, which it is equal to when $\rho=22$)\footnote{Now that 
Artin's conjecture has been proved for $p\geq 5$ 
one might work as well with ${\NS}\Tensor\ZZ_p$}. We let $N$ be the flat
cohomology group $H^2(X,\ZZ_p(1))$ and consider $N\Tensor {\bf k}$ 
with $F$ acting as
${\rm id}\tensor F$.  We then have de Rham Chern class map
$c_1: {N}\to H^2_{dR}(X/{\bf k})$ 
and we shall also write $c_1$ for the ${\bf k}$-linear
extension $N\Tensor{\bf k} \to H^2_{dR}(X/{\bf k})$ of 
$c_1$. The kernel of this map is called characteristic subspace and
plays the central role, see \cite{ogus79::super+k3}. 
We write the kernel of it on the form $F^*K$ for
some sub-vector space $K \subseteq N\Tensor{\bf k}$.  We let $\tilde K$ be the
inverse image of $K$ in $N\Tensor{\bf W}$.  We can consider ${\cH}$ as a
${\bf W}$-submodule of $N\Tensor Q$, where $Q$ is the fraction field of ${\bf W}$.
Then, by definition and the fact that $p\dual{N}\subseteq N$, with
$\dual{N}$ the dual of $N$ with respect to the
intersection pairing, we have that ${\cH}=p^{-1}\sigma^*\tilde{K}$. Furthermore,
as $F=p\tensor \sigma$ on $N\Tensor{\bf W}$ it is clear that
$$
{\cH}_1 = F^{-1}p {\cH}=p^{-1}(\tilde{K}\cap\sigma^*\tilde{K})
\quad {\rm and}\quad 
{\cH}_0 =F^{-1}p^2K=\tilde{K}
$$
and they map to the Hodge filtration of $H$. On the other hand,
$$
{\cH}^{c}_0 =F({\cH}_2)=\sigma^{*2}(\tilde K) \quad {\rm and} \quad
{\cH}^c_1=p^{-1}F({\cH}_1)=p^{-1}(\sigma^*\tilde{K}\cap\sigma^{*2}\tilde{K})
$$
which then map to the conjugate filtration.
Starting our procedure for constructing
the canonical filtration we see that it stops immediately when
$H_0=H^c_0$, but this is the case precisely when the Artin
invariant $\sigma_0$ equals $1$.
If not, we add the image $\overline{E}_2$ of
$H^c_0$ in $H/H_0$ to the Hodge filtration, whose inverse
image in ${\cH}$ then is $\tilde
U_2:=\tilde{K}+\sigma^{*}\tilde{K}+\sigma^{*2}\tilde{K}$.
The next step is to transfer $\overline{E}_2$ via the Cartier
isomorphism to get an addition, $V_2$,
to the conjugate filtration. As the Cartier isomorphism between the ``middle
parts'' of the Hodge and conjugate filtrations is implemented by
$p^{-1}F=\sigma^*$ we get that the inverse image of $V_2$ in $H$ is given by
$$
\sigma^*\tilde{K}+\sigma^{*2}\tilde{K}+\sigma^{*3}\tilde{K}.
$$
The process stops at that stage precisely when $\sigma^{*3}\tilde{K}\subseteq
\tilde{K}+\sigma^*\tilde{K}+\sigma^{*2}\tilde{K}$ which in turn is equivalent to
$\tilde{K}+\sigma^{*}\tilde{K}+\sigma^{*2}\tilde{K}$ being stable under
$\sigma^*$. However, that in turn is equivalent to
$\tilde{K}+\sigma^{*}\tilde{K}+\sigma^{*2}\tilde{K}=pN^*\Tensor{\bf W}({\bf k})$ (by
\cite[3.12.3]{ogus79::super+k3}) and thus to $\sigma_0=2$ (as we have put
ourselves in the case when $\sigma_0>1$). If not, the process continues, forcing
us to add the image of
$\sigma^*\tilde{K}+\sigma^{*2}\tilde{K}+\sigma^{*3}\tilde{K}$ to the Hodge
filtration. If we continue in this way it is clear that we will stop at
$\tilde{K}+\sigma^{*}\tilde{K}+\cdots+\sigma^{*\sigma_0}\tilde{K}$ which equals
$pN^*\Tensor{\bf W}({\bf k})$. In this way we get an extension of the Hodge
filtration which ends at $R\Tensor{\bf k}$, where $R$ is the radical of
$\overline{N}:=N\Tensor{\FF}_p$. Its annihilator is $N\Tensor{\bf k}$ and hence we
get that the ``middle subquotient'' of the canonical filtration is canonically
isomorphic to $\overline{N}/R\Tensor{\bf k}$ with the natural quadratic structure
and the map induced by the Cartier isomorphism having $\overline{N}/R$ as its
fixed points. Hence extending the canonical filtration to a final one amounts to
finding a complete self-dual flag in $\overline{N}/R$. We also get from
\cite[3.4]{ogus79::super+k3} and the fact that the discriminant of $N(X)$ is
$-p^{2\sigma_0}$ that the quadratic form on $\overline{N}/R$ is non-split.
Also in this case it is evident by inspection
from the filtration thus obtained that it is of type $w_{n-1-\sigma_0}$
or $w_{n-1-\sigma_0}^{\prime}$. In the case where $\sigma_0=n/2$ we find
 $w_{n-1-\sigma_0}$. In the general case we see that the $F$-zip is the sum of
 an $F$-zip of dimension $2\sigma_0$ and a Tate $F$-zip which is isomorphic to
 the middle part $F$-zip. From the multiplicativity of the Hodge discriminant
 and the case $n=2\sigma_0$ we conclude.
\end{proof}
\end{theorem}

\end{section}

\begin{section}{Strata on the Flag Space}
Here we shall define strata on the flag space of orthogonal or symplectic
flags on the primitive part of the second de Rham cohomology of a family 
of polarized K3 surfaces. Hence we assume that we have a family
  $f\colon X\to S$ of $N$-marked K3 surfaces (where $S$ may be an algebraic
  stack). The primitive cohomology forms a vector bundle ${\mathcal H}$ 
over $S$ of rank $n$ with an orthogonal structure given by the
intersection form. It is provided with two orthogonal partial flags: the
Hodge flag and the conjugate flag, thus giving a $F$-zip. 

If we choose an orthogonal flag refining
the Hodge filtration $C^{\bullet}$ 
we obtain by using the Cartier isomorphism $\varphi_i: F^*(C^i/C^{i+1}) \cong
D_{i+1}/D_i$ a second
flag refining the conjugate filtration $D_{\bullet}$. 

We let $\cF_n$ be the space of complete 
orthogonal flags on ${\mathcal H}$ refining the Hodge filtration. 
It admits a canonical projection $\cF_n \to S$. Since a flag refining
the Hodge filtration automatically defines a second flag (refining the
conjugate filtration) we can measure the relative position of these flags
and thus define strata on $S$. We refer to \cite{FP98} for background.

We can formulate this in the following way, cf.\ \cite{ekedahl10::cycle+class+e+o+strat}.
Let $G$ be a semi-simple group
and $B$ a Borel subgroup and $G/B$-bundle $T \to Y$ over some scheme $Y$ 
with $G$ as structure group. Suppose that we have two sections $t_i: Y\to
T$ of $T$ with $i=1,2$. If $w$ is an element of the Weyl group of $G$ we define 
a locally closed subscheme ${\Uo_w}$ (resp.\ ${\Uc_w}$) of $Y$ by
choosing locally (possibly in the \'etale topology) 
a trivialization of $T$ for which $t_1$ is a constant section.
Then $t_2$ corresponds to a map $Y \to G/B$ and we define ${\Uo_w}$ (resp.\ ${\Uc_w}$) to be the inverse image of the $B$-orbit $BwB$ (resp.\ of its closure
in $G/B$). This does not depend on the trivialization taken since the 
difference corresponds to a map $Y \to B$ and the cycles $BwB$ and its
closure are $B$-invariant. Therefore this defines global subschemes 
 ${\Uo_w}$ (resp.\ ${\Uc_w}$) of $Y$. If $t_1$ and $t_2$ have the property
that $Y={\Uc_w}$ then we say that $t_1$ and $t_2$ are in relative position $w$.
We thus find our strata $\Uo_w$ and $\Uc_w$ associated to an element of our
Weyl group. Note that a priori it is not clear that the closure of
$\Uo_w$ equals $\Uc_w$, but this will hold (see below). 

We shall apply this to the situation that $\cF_n$ is the space of orthogonal
flags refining the conjugate filtration on 
${\mathcal H}$ for the family $f:X \to S$. 
\begin{definition}
On the space ${\mathcal F}_n$ of orthogonal flags on ${\mathcal H}$
we define for every element $w$ in our Weyl group 
the locally closed subschemes ${\Uo_w}$ and ${\Uc_w}$ of ${\mathcal F}_n$ 
associated to the flag refining the conjugate filtration and the induced
flag on the Hodge filtration as the subschemes that measure the relative
position of these two orthogonal flags as defined above. 
For final elements $w$ in our Weyl group
we have an orthogonal flag refining the canonical filtration 
and then define the stratum ${\cV}_w$ on $S$ as the locally
closed subset of $S$ of points for which the canonical type of the K3 surface
is equal to the canonical type of $w$. By Thm.\ 
\ref{Final/canonical filtrations type} these strata ${\cV}_w$ 
belong to the height and
Artin invariant. 
\end{definition}

It might seem that working on the flag space ${\mathcal F}_n$ rather than on $S$
is a detour,
but in the next section we shall see that it helps understanding the strata
on the moduli.

\end{section}

\begin{section}{The local structure of strata}\label{localstructure}

The reason for working on the flag space over our moduli space (of lattice polarized K3 surfaces) is that the
  strata are much better behaved than on the moduli space. In fact, up to
  infinitesimal order $<p$ the strata look like usual Schubert strata.  This
  idea of \cite{ekedahl10::cycle+class+e+o+strat} and the methods employed there
  can be transferred to our situation. Hence we assume that we have a family
  $f\colon X\to S$ of $N$-marked K3 surfaces (where $S$ may be an algebraic
  stack). We shall also need to assume a versality condition: For any geometric
  point $s$ of $S$ contraction of forms by vector fields induces a map
  $H^1(X_s,T_{X_s})\to {\rm Hom}(H^0(X_s,\Omega^2_{X_s}),H^1(X_s,\Omega^1_{X_s}))$ we
  can then compose this with the map induced by the projection on the second
  factor of the decomposition $H^1(X_s,\Omega^1_{X_s})=N\Tensor{\bf k}
\perp P$ and
  then further compose the resulting map with the Kodaira-Spencer map $T_sS\to
  H^1(X_s,T^1_{X_s})$. The required versality condition is that $S$ be smooth at
  all $s$ and that the composed map $T_sS\to {\rm Hom}(H^0(X_s,\Omega^2_{X_s}),P)$ be
  surjective.

  The space $\cF_n$ together with the $\Uo_w$ is a stratified space.  The space
  $\Fl_n$ of complete self-dual flags on an orthogonal space is also a
  stratified space with the stratification given by the Schubert cells. The idea
  is that our flag space at a point can be identified up to the $(p-1)$st
  neighborhood with the flag space at an appropriate point. Moreover, under this
  correspondence the strata on $\cF_n$ correspond precisely to the Schubert
  strata on ${\Fl}_n$. This enables us to transplant the detailed knowledge
  about Schubert strata up to order $p$ to our situation. More precisely, if $R$
  is a local ring with maximal ideal $m$ defining an affine scheme $S$ then the
  height-$1$ hull of $R$ (or $S$) is given by $R/m^{(p)}$, with $m^{(p)}$
  generated by the $p$'th powers of elements of $m$.  We call two local rings
  {\sl height $1$-isomorphic} if their respective height $1$-hulls are
  isomorphic.
\begin{theorem}
  Let $k$ be a perfect field of positive characteristic $p$.  For each
  $k$-point $x$ of $\cF_n$ there is a $k$-point $y$ of $\Fl_n$ such that the
  height $1$-neighbourhood of $x$ is isomorphic to the height $1$-neighbourhood
  of $y$ times a smooth space by an isomorphism respecting stratifications.
\begin{proof}
We consider the de Rham cohomology ${H}$ together with the
Gauss-Manin connection on the height-$1$ neighborhood $Y$ of $x$.
We can trivialize ${H}$ plus its Gauss-Manin connection on $Y$
since the ideal of $x$ has a divided power structure for which
divided powers of degree $\geq p$ are zero. This implies that the
orthogonal flags $E_{\bullet}$ and $G_{\bullet}$ on ${H}$ are
horizontal. We thus get a map from $Y$ to the space of orthogonal
flags on a standard orthogonal space, that is, an isomorphism
from $Y$ to a height-$1$ neighborhood on $\Fl_n$. It is not difficult
to see that it preserves strata.
\end{proof}
\end{theorem}
The following theorem is the main consequence of this. Denote the base
space of $\cF_n$ by $\cK_n$.

\begin{theorem} The strata $\Uo_w$ possess the following properties:
\begin{enumerate}
\item[i)] Each stratum $\Uo_w$ is smooth of dimension $\ell(w)$.
\item[ii)] The closed stratum $\Uc_w$ is reduced, Cohen-Macaulay and normal of dimension
$\ell(w)$ and is the closure of $\Uo_w$ for all $w$ in the Weyl group.
\item[iii)] If $w$ is final then the restriction to $\Uo_w$ 
of the projection $\cF_n \to \cK_n$ to $\Uo_w$ is a
finite surjective \'etale covering from $\Uo_w$ to $\sV_w$ of degree
equal to the number of final filtrations on a canonical filtration
of type $w$.
\end{enumerate}
\begin{proof}
  This theorem follows from Thm 10.1 in exactly the same fashion as Corollary
  8.4 in Section 8.2 of \cite{ekedahl10::cycle+class+e+o+strat} follows from
  Theorem 8.1 there.
\end{proof}
\end{theorem}
In view of our calculations of the cycle classes we need to know the degrees
of the canonical projections $\pi_w: \Uo_w \to \sV_w$ of our strata for
final or twisted final elements. This
degree is expressed as the number of ways we can put a final filtration
on a canonical one. We now calculate these degrees.

\begin{lemma}\label{DegreelemmaB}
Let $w\in \Weyl{B}{m}$ be a final element and let $\pi_w : \Uo_w \to \sV_w$
be the restriction of the projection from $\cF_n\to \cK_n$ with $n=2m+1$.
\begin{enumerate}
\item[i)] For $1\leq k \leq m-1$ we have
$\deg(\pi_{w_k})/\deg(\pi_{w_{k+1}})=p^{2m-2k-1}+p^{2m-2k-2}+\ldots + 1$.
\item[ii)] Similarly, we have
$\deg(\pi_{w_{m+k+1}})/\deg(\pi_{w_{m+k}})=p^{2k-1}+p^{2k-2}+\cdots+1$.
\end{enumerate}
\begin{proof}
Note that $\deg(\pi_{w})$ is the number of final filtrations
on a given canonical filtration of type $w$.
For case i) we look at the number of lines in a linear space of dimension
$2m-2k$ over ${\FF}_p$, i.e.\ the number of points in projective space of
dimension $2m-2k$. For case ii) we look at the number of isotropic lines
in an orthogonal space of dimension $2k+1$, i.e.\ the degree is the
number of points on a quadric of dimension $2k-1$.
\end{proof}
\end{lemma}
\begin{lemma}\label{DegreelemmaD}
Let $w\in \Weyl{D}{m}$ be a final element and let $\pi_w : \Uo_w \to \sV_w$
be the restriction of the projection from $\cF_n\to \cK_n$ with $n=2m$.
\begin{enumerate}
\item[i)] For $1\leq k \leq m-1$ we have
$\deg(\pi_{w_k})/\deg(\pi_{w_{k+1}})=-p^{m-k-1}+\sum_{j=0}^{2m-2k-2} p^j$.

\item[ii)] We have $\deg(\pi_{w_{m-1}})=
\deg(\pi_{w_{m}})=(\pi_{w_{m+1}})=(\pi_{w_{m+2}})=1$.

\item[iii)] Similarly, for $2\leq k \leq m-1$ we have
$\deg(\pi_{w_{m+k+1}})/\deg(\pi_{w_{m+k}})=p^{k-1}+\sum_{j=0}^{2k-2} p^j$.
\end{enumerate}
\begin{proof}
The proof is the same as for Lemma \ref{DegreelemmaB} except that the counts of
isotropic lines are different. It also  depends on whether the form is split or
not but that is provided by Theorem \ref{Final/canonical filtrations type}.
\end{proof}
\end{lemma}

\begin{lemma}\label{DegreelemmaDtwisted}
Let $w=w^{\prime}s_m^{\prime} \in \Weylp{D}{m}$ be a twisted
final element and  let $\pi_w : \Uo_w \to \sV_w$ be the
restriction of the projection from $\cF_n\to \cK_n$ with $n=2m$.
\begin{enumerate}
\item[i)] For $1\leq k \leq m-1$ we have
$\deg(\pi_{w_k})/\deg(\pi_{w_{k+1}})=p^{m-k-1}+\sum_{j=0}^{2m-2k-2} p^j$.

\item[ii)] We have $\deg(\pi_{w_{m-1}})=
\deg(\pi_{w_{m}})=(\pi_{w_{m+1}})=(\pi_{w_{m+2}})=1$.

\item[iii)] Similarly, for $2\leq k \leq m-1$ we have
$\deg(\pi_{w_{m+k+1}})/\deg(\pi_{w_{m+k}})=-p^{k-1}+\sum_{j=0}^{2k-2} p^j$.
\end{enumerate}
\begin{proof}
  Again the proof is the same as for Lemma \ref{DegreelemmaB} with needed extra
  information provided by Theorem \ref{Final/canonical filtrations
    type}.
\end{proof}
\end{lemma}

\end{section}

\begin{section}{Shuffles}

In this section we shall discuss an analogue and generalization of shuffles
introduced in 
\cite{ekedahl10::cycle+class+e+o+strat}. 
They play a role
in describing maps between the strata $\Uo_w$ on the flag space and
describe inseparable maps between strata.
They are a key instrument for deciding whether the push forward of
the corresponding cycle class will vanish. In fact, we saw that for a
final stratum $\Uo_w$ 
the projection to $\sV_w$ is finite \'etale. It will turn out
that for a non-final stratum the projection is lower-dimensional or
factors through an inseparable map to a final stratum. These inseparable
maps are described by shuffle maps as in \cite{ekedahl10::cycle+class+e+o+strat}. We analyze the situation in detail.

\smallskip
We start with an elementary lemma on the length of elements in our Weyl groups.
A general reference is \cite{BB05}.

\begin{lemma}\label{length} The length satisfies the following properties.
\begin{enumerate}
\item[i)] Let $w$ be an element of a Coxeter group and suppose that
$\ell(ws_i)=\ell(w)-1$. Then either $\ell(s_iws_i)=\ell(w)$ or
$\ell(s_iws_i)=\ell(w)-2$.

\item[ii)] For $w$ an element of $\Weyl{B}{m}$, $\Weyl{C}{m}$ or $\Weyl{D}{m}$
and $1 \leq i < m$ we have that $\ell(ws_i)=\ell(w)-1$
precisely when $w(i+1)<w(i)$ and then
$\ell(s_iws_i)=\ell(w)-2$ precisely when $w^{-1}(i+1)<w^{-1}(i)$.

\item[iii)] For $w$ an element of $\Weyl{B}{m}$ we have $\ell(ws_m)=\ell(w)-1$
precisely when $w(m)>w(m+1)$ and then $\ell(s_mws_m)=\ell(w)-2$
precisely when $w^{-1}(m+2) < m$.

\item[iv)] For $w$ an element of $\Weyl{C}{m}$ we have $\ell(ws_m)=\ell(w)-1$
precisely when $w(m)>w(m+1)$ and then $\ell(s_mws_m)=\ell(w)-2$
precisely when $w^{-1}(m+1)<m$.

\item[v)] For $w \in \Weyl{D}{m}$ we have $\ell(ws_m)=\ell(w)-1$ precisely when
$w(m-1)>w(m+1)$ and $w(m)>w(m+2)$ and then
$\ell(s_mws_m)=\ell(w)-2$ precisely when $w^{-1}(m+1)>w^{-1}(m-1)$
and $w^{-1}(m+2)> w^{-1}(m)$.
\end{enumerate}
\begin{proof}
Easy. 
\end{proof}
\end{lemma}
We now define the shuffle maps.
Assume that the dimension of our orthogonal or symplectic space is $n$ 
and that $W_n$ denotes the Weyl group in question. 
Assume that $w \in {W}_n$ and that
$\ell(ws_i)=\ell(w)-1$ for some $1<i\leq m$.
This means that for the universal flags $\EE_\bullet$ and $\GG_\bullet$
on $\Uo_w$ the image of $G_{w(i+1)}\cap \EE_{i+1}$ in
$\EE_{i+1}/\EE_{i-1}$ is a line bundle. We define a new self-dual flag
$\EE'_\bullet$ on $\Uo_w$ by the condition that
$\EE'_j=\EE_j$ for $i \neq j \leq n$ and $\EE'_i/\EE_{i-1}$ be equal
to the image of $\GG_{w(i+1)}\cap \EE_{i+1}$.
This then gives a map
$$
\sigma_{w,i}: \Uo_w \to \cF_n, \qquad (\EE_{\bullet},\GG_{\bullet})
\mapsto (\EE^{\prime}_{\bullet},\GG_{\bullet})
$$
which we shall call the {\sl $i$'th elementary shuffle map} for $w$.
We say that the shuffle map is {\sl unambiguous} if there is a
$v \in {W}_n$ such that the image of $\sigma_{w,i}$ lies in $\Uo_v$.
Please note the condition $i>1$ which ensures that the first and
last step of the Hodge filtration are left unchanged.
\begin{proposition}\label{ambiguous}
The shuffle map $\sigma_{w,i}$ satisfies the following properties.
\begin{enumerate}
\item[i)] The element $\sigma_{w,i}$ for $w \in {W}_n$ is unambiguous precisely
when $\ell(s_iws_i)=\ell(w)$. In that case the image of $\sigma_{w,i}$ is equal
to $\Uo_{s_iws_i}$ and $\sigma_{w,i}$ is finite and purely inseparable of degree
$p$.
\item[ii)] If $\ell(s_iws_i)=\ell(w)-2$ then $\sigma_{w,i}$ maps onto
$\Uo_{ws_i}\cup\Uo_{s_iws_i}$. In particular it is not generically finite.
\end{enumerate}
\begin{proof}
Assume first that $\ell(s_iws_i)=\ell(w)$. We may locally (in the \'etale
topology) choose a basis adapted to the two flags, i.e., an orthonormal basis
$e_1,\dots,e_n$ of $M$ (on $\Uo_w$) such that $\EE_j$ is spanned by
$e_1,\dots,e_j$ and $\GG_{w^{-1}(j)}$ is spanned by
$e_{w^{-1}(1)},\dots,e_{w^{-1}(j)}$. We then have that $\EE'_i$ is spanned by
$e_1,\dots,\hat{e}_i,e_{i+1}$. We may further assume that
$C^{-1}e_j=e_{w^{-1}(j)}\bmod \GG_{j-1}$ for all $1<j<n$. Put $k := w^{-1}(i)$
and $\ell:=w^{-1}(i+1)$, we then have, by the assumption and Lemma \ref{length},
that $k<\ell$. There is a $\lambda$ such that $C^{-1}e_{i+1}=e_\ell+\lambda
e_k\bmod \GG_{i-1}$. We now put for $j\leq m$
\begin{displaymath}
e'_j =
\begin{cases}
e_j&\text{if $j\neq i,i+1,\ell$,}\\
e_{i+1}&\text{if $j=i$,}\\
e_i    &\text{if $j=i+1$, and}\\
e_\ell+\lambda e_k&\text{if $j=\ell$.}
\end{cases}
\end{displaymath}
By the assumption $k<\ell$ we get that $\EE'_j$ is spanned by $e'_1,\dots,e'_j$
and we may extend it (uniquely) to an adapted basis for $\EE'_\bullet$ and
$\GG'_\bullet$. This makes it clear that $(\EE'_\bullet,\GG'_\bullet)$ is of
type $s_iws_i$.
Consider conversely the universal flag $\EE_\bullet$ on $\Uo_v$, where
$v=s_iws_i$ and put again $k := v^{-1}(i)$ and $\ell:=v^{-1}(i+1)$,
where this time $k>\ell$ and $v^{-1}(i)<v^{-1}(i+1)$.
We choose as before an adapted basis and  let $\EE'_i$ be spanned by
$e_1,\dots,e_{i+1}+\rho e_i$ and then $\GG'_i$
is spanned by $e_\ell+(\rho^p+\lambda)e_k$ and $\GG_{i-1}$.
It is then easy to see that $(\EE'_\bullet,\GG'_\bullet)$ is of type
$w$ precisely when $\rho^p+\lambda=0$ which gives i).

We now assume that $\ell(s_iws_i)=\ell(w)-2$. The setup is then the same as
before except that we now have $k>\ell$. This means that
$(\EE'_\bullet,\GG'_\bullet)$ will be of type $ws_i$
if $\lambda \neq 0$ and of type $s_iws_i$ if $\lambda=0$.
The converse is similar to the converse of  i)
(with the difference that the new flag pair will be of type $w$
for all choices of $\rho$).
\end{proof}
\end{proposition}
\begin{definition}
A sequence of elements $w_1,\dots,w_r$ of ${W}_n$ is said to be a
{\sl shuffle sequence} if for each $1\leq k<r$ there is an $1<i_k\leq m$ such
that $\ell(w_ks_{i_k})=\ell(w_k)-1$ and $w_{k+1}=s_{i_k}w_ks_{i_k}$
with $\ell(s_{i_k}w_ks_{i_k})=\ell(w_k)$. It
is said be {\sl ambiguous} if there is an $i_r$ such that $\ell(w_r
s_{i_r})=\ell(w_r)-1$ and $\ell(s_{i_r}w_ks_{i_r})=\ell(w_r)-2$,
{\sl final} if $w_k$ is final, and {\sl cyclic}
if there are $1\leq j <k\leq r$ such that $w_j=w_k$.
\end{definition}
We shall now show that starting with a non-final element we always can
find a shuffle sequence. Doing this we see that we 
end up at a final stratum or
we find that the image of our stratum under projection is lower-dimensional. 
Recall that $\cK_n$ denotes the base space of our flag space $\cF_n$.
\begin{proposition}\label{shufflekinds}
Shuffle sequences exist:
\begin{enumerate}
\item[i)] For every $w \in {W}_n$ there exists a shuffle 
sequence starting with $w$
and which is either ambiguous, final, or cyclic.

\item[ii)] If there is an ambiguous or cyclic shuffle sequence 
starting with $w \in {W}_n$, then the restriction of the 
projection map $\cF_n \to \cK_n$ to $\Uo_w$ is
not generically finite.
\item[iii)] Given a final shuffle sequence $w_1,\dots,w_r$, the restriction of
the projection map $\cF_n \to \cK_n$ to $\Uo_{w_1}$ is the composite of a finite
flat purely inseparable map of degree $p^{r-1}$ and the finite \'etale map
$\Uo_{w_r} \to \Vo_{w_r}$.
\end{enumerate}
\begin{proof}
Let $w=w_1,\dots,w_r$ be a shuffle sequence which is maximal for not being
ambiguous, final, or cyclic. In particular $w_r$ is not final and therefore there
is an $1<i_r\leq m$ such that $\ell(w_rs_{i_r})=\ell(w_r)-1$. If
$\ell(s_{i_r}w_rs_{i_r})=\ell(w_r)$, then by the maximality we must have either that
$w_{i_{r+1}}:=s_{i_r}w_rs_{i_r}$ is final or appears in the sequence so that we
get a final or cyclic sequence by adding $w_{i_{r+1}}$. If
$\ell(s_{i_r}w_rs_{i_r})=\ell(w_r)-2$ we instead get an ambiguous sequence thus
proving part i).

If there is an ambiguous sequence $w=w_1,\dots,w_r$, then the projection map
$\Uo_w$ factors by Proposition \ref{ambiguous} as $\Uo_w \to \Uo_{w_rs_{i_r}}\cup
\Uo_{s_{i_r}w_rs_{i_r}} \to \cK_n$ and as
$\ell(w_rs_{i_r})< \ell(w)$ and $\ell(s_{i_r}w_rs_{i_r})<\ell(w)$
and hence $\Uo_w$ has an
image of dimension smaller than that of $\Uo_w$.
On the other hand if there is a cyclic sequence,
then the projection factors through an infinite sequence of
$\sigma_{v,j}$'s and as each of them is of degree $>1$ we get that image has
lower dimension. This proves part ii).

Finally, assume that we have a final sequence $w=w_1,\dots,w_r$. Then the
projection factors as the composite
$\sigma_{w_{r-1},i_{r-1}}\circ\cdots\circ\sigma_{w_{1},i_{1}}$ and the
projection $\Uo_{s_{i_r}w_rs_{i_r}} \to \cK_n$. The latter is an \'etale cover
of $\Vo_{w_r}$ and the first is finite purely inseparable of degree $p$.
\end{proof}
\end{proposition}
We shall call an ambiguous or cyclic shuffle a {\sl degenerate
shuffle}. Proposition \ref{shufflekinds} implies
that either all shuffles of an element $w \in {W}_n$ are
degenerate or they are all final. In the first case, the
projection map restricted to $\Uc_w$ is not generically finite on
each of its irreducible components and in particular the image of
$[\Uc_w]$ is zero. In the second the class of the push forward
is non-zero and equal $p^\ell [\Uc_\nu]$, where $\ell$ is the
length of a final shuffle of $w$ to the final element $\nu$.
\end{section}
\begin{section}{Final elements}
In order to calculate cycle classes of our strata we shall apply
a Pieri formula which gives an expression of the intersection
product of a class of a stratum with a first Chern class in terms
of cycle classes of strata of dimension one less. For this we need
a precise description of the colength one elements in the Weyl group
below a given final (or twisted final) element and then determine
whether these elements are degenerate or of shuffle type.
In this auxiliary and rather technical section we describe the elements
involved. After treating the case of $\Weyl{B}{m}$ in detail we deal
with the other cases more succinctly. 

\begin{subsection}{Final elements in $\Weyl{B}{m}$}
We begin by factoring the final elements in the Weyl group $\Weyl{B}{m}$ as a
product of simple reflections.
\begin{lemma}\label{final factorisations}
The products $w_k=s_ks_{k+1}\cdots s_ms_{m-1}\cdots s_1$
with $1\leq k \leq m-1$ and $w_{2m-k}=s_ks_{k-1}\cdots s_1$
with $m\geq k \geq 0$ are reduced expressions for the
$2m$ final elements of $\Weyl{B}{m}$. We have $w_1=w_{\emptyset}$ and
$w_{2m}=1$.
\begin{proof}
Easily verified.
\end{proof}
\end{lemma}
Note that the final elements are linearly ordered by the their length.
We now determine the elements of colength $1$ below a final
element in the Bruhat order.
\begin{proposition}
The elements of colength $1$ below a final
element of $W^B_m$  in the Bruhat order are as follows:
\begin{enumerate}
\item[i)] The elements in $W^B_m$ of colength $1$ below the final
element $w=s_k\cdots s_m\cdots s_1$ with $k <m$ are
$s_k\cdots \hat{s}_i\cdots s_m\cdots s_1$ for $i=k,\ldots,m-1$,
and $s_k\cdots s_m\cdots \hat{s}_i \cdots s_1$ for $i=m-1,m-2,\ldots,1$.
They are obtained from the final element by multiplying $w$
to the right by the element $s_\alpha$, where $\alpha$ is the root
$\epsilon_1+\epsilon_{k+1}$,\dots,$\epsilon_1+\epsilon_m$,
$\epsilon_1-\epsilon_m$,\dots,$\epsilon_1-\epsilon_2$ respectively.

\item[ii)] The elements of colength one below $w=s_ms_{m-1}\cdots s_1$
are the elements $s_m\cdots \hat{s}_i\cdots s_1$ for $i=m,\ldots,1$.
They are obtained by
multiplying $w$ from the right by $s_{\alpha}$ where $\alpha$ is the root
$\epsilon_1$, $\epsilon_1-\epsilon_m$, $\epsilon_1-\epsilon_{m-1}, \ldots,
\epsilon_1-\epsilon_2$ respectively.

\item[iii)] The elements of colength one below $w=s_{k}s_{k-1}\cdots s_1$
are the elements $s_k\cdots
\hat{s}_i\cdots s_1$ for $i=k,\ldots,1$.
They are obtained from the final element by multiplying
$w$ to the right by the element $s_\alpha$, where $\alpha$ is the root
$\epsilon_1-\epsilon_{k+1}$,\dots,$\epsilon_1-\epsilon_2$ respectively.
\end{enumerate}
\begin{proof}
We know that the elements of colength $1$ below an element are obtained by
considering a reduced expression for the element, taking the elements obtained
by removing one element from the expression, and then keeping the elements of
colength $1$. Lemma \ref{final factorisations} provides a reduced expression. In
case part {i)}, among the elements obtained by removing one simple
reflection from the reduced expression clearly the one obtained by removing
$s_m$ (when present) is not of colength $1$ and the others are easily shown to
be. Finally, if the element has the factorization $w's_iw''$ and the colength
$1$ element has the factorization $w'w''$ then it is obtained by multiplying by
$(w^{\prime\prime})^{-1}s_iw^{\prime\prime}$ to the right , i.e., by $s_\alpha$,
where $\alpha = (w'')^{-1}(\alpha_i)$. From this the rest follows by a simple
calculation. The cases part {ii)} and  part {iii)} are similar.

\end{proof}
\end{proposition}
We shall now consider the elements of colength $1$ below a final element and
determine if they have degenerate or final shuffle type.
\begin{proposition}
The elements of colength $1$ below a final element satisfy the following.
\begin{enumerate}
\item[i)] The element $s_ks_{k+1}\dots s_m\cdots\hat{s}_i\cdots s_1$ with $1\leq
i<k<m$ is degenerate.

\item[ii)] The element $s_ks_{k+1}\cdots s_m\cdots\hat{s}_i\cdots s_1$
with $m>i\geq k$ has an
elementary shuffle to the element $s_ks_{k+1}\cdots s_m\cdots\hat{s}_{i+1}\cdots
s_1$ if $i<m-1$ and to $s_ks_{k+1}\cdots \hat{s}_{m-1}s_m\cdots s_1$ if $i=m-1$.

\item[iii)] For $m>i>k$ the element $s_k\cdots\hat{s}_i\cdots s_m\cdots s_1$
has an elementary shuffle to
the element $s_k\cdots \hat{s}_{i-1}\cdots s_m\cdots s_1$.

\item[iv)] The element $s_k\cdots\hat{s}_i\cdots s_1$ is degenerate if $i<k
\leq m$.
\end{enumerate}
\begin{proof}
Starting with part {i)} we note that the elements $s_j$ with
$k\leq j \leq i$ commute with the $s_l$ with $i-1 \geq l \geq 1$.
This means that $s_k\cdots s_m\cdots\hat{s}_i\cdots s_1 =
s_k\cdots s_m\cdots s_{i+2}s_{i-1}\cdots s_1s_{i+1}$
and this implies that we can perform an $i+1$'st shuffle giving the element
$s_{i+1}s_k\cdots s_m\cdots s_{i+2}s_{i-1}\cdots s_1$.
If $i+1=k$ this element has shorter length, while if
$i+1=k-1$ we get $s_{k-1}\cdots s_m \cdots s_{i+2}s_{i-1}\cdots s_1$
and then move $s_{i+2}$ to the right and perform a shuffle with it.
We thus arrive at the element
$s_ks_{k-1}s_k\cdots s_m \cdots s_{i+3}s_{i-1}\cdots s_1$,
and by applying the braid relation we get the element
$s_{k-1}s_{k}s_{k-1}s_{k+1}\cdots s_m \cdots s_{i+3}s_{i-1}\cdots s_1$
and by moving $s_{k-1}$, the third factor, this is seen to equal
$s_{k-1}s_{k}\cdots s_m \cdots s_{i+3}s_{i+1}s_{i-1}\cdots s_1$,
and by moving $s_{i+1}=s_{k-1}$ to the right, then
performing a shuffle by $s_{i+1}$
we get an element of shorter length. If however, $i+1<k-1$ we get
$s_k\cdots s_m \cdots s_{i+3}s_{i+1}s_{i+2}s_{i-1}\cdots s_1$
and then we can perform a shuffle by $s_{i+2}$. Continuing
in this way this leads to a shorter element.

We continue with part {ii)} and assume that $i+1\neq m$. Then
$s_ks_{k+1}\cdots s_m\dots s_{i+1}s_{i-1}\cdots s_1$ equals
$ s_ks_{k+1}\cdots s_m\cdots s_{i+2}
s_{i-1}\cdots s_1s_{i+1}$ as the involved simple transpositions commute. This
means that we may perform an elementary shuffle to get the element
$s_{i+1}s_ks_{k+1}\cdots s_m\cdots s_{i+2}s_{i-1}\cdots s_1$
which in turn is equal to
$s_ks_{k+1}\cdots s_{i+1}s_is_{i+1}\cdots s_m\cdots s_{i+2}s_{i-1}\cdots s_1$.
Using the braid rule we get that this element equals
$s_ks_{k+1}\cdots s_{i}s_{i+1}s_{i}\cdots
s_m\cdots s_{i+2}s_{i-1}\cdots s_1$ and we observe that this
 in turn equals $s_ks_{k+1}\cdots s_m \cdots s_{i+2}s_is_{i-1}\cdots s_1$
and this is
$s_{k}\cdots s_m \hat{s}_{i+1}\cdots s_1$. By Lemma \ref{final
factorisations} this element is reduced so that we have performed an unambiguous
shuffle to the claimed element. If instead $i+1=m$ we have that
$s_ks_{k+1}\cdots s_m s_{m-2}\cdots s_1$ is equal to $s_ks_{k+1}\cdots
s_{m-1}s_{m-2}\cdots s_1s_m$
which an elementary shuffle turns into $s_ms_ks_{k+1}\cdots
s_{m-1}s_{m-2}\cdots s_1$ which on its turn
equals the element  $s_ks_{k+1}\cdots
s_{m-2}s_ms_{m-1}s_{m-2}\cdots s_1$, a reduced expression of the right
element.

For part {iii)} note that
 $s_k\cdots s_{i-1}s_{i+1}\dots s_m\cdots s_1$
equals $s_k\cdots s_{i-2}s_{i+1}\cdots s_m\cdots s_{i-1}s_{i}s_{i-1}\cdots s_1$
by moving $s_{i-1}$ to the right, and
by the braid rule equals $s_k\cdots s_{i-2}s_{i+1}\cdots s_m\cdots
s_{i}s_{i-1}s_{i}\cdots s_1$ in which the last
$s_i$ migrates to the right to give
 $s_k\cdots s_{i+1}\cdots s_m\cdots
s_{i}s_{i-1}\cdots s_1s_i$. Performing an elementary shuffle leads to
$s_is_k\cdots s_{i+1}\cdots s_m\cdots s_{i}s_{i-1}\cdots s_1$ which is equal to
the element
$s_k\cdots s_is_{i+1}\cdots s_m\cdots s_{i}s_{i-1}\cdots s_1$. This is, still by the
lemma, a reduced expression of the desired element.
The proof of part {iv)} is analogous to that of part {i)}.
\end{proof}
\end{proposition}
\begin{corollary}
The non-degenerate elements of colength $1$ below the final ones are as follows.
\begin{enumerate}
\item[i)] For a final element $w=s_k\ldots s_m\ldots s_1$ with $k<m$
the only non-degenerate elements of colength $1$ below $w$ are
$ws_{\alpha}$ with $\alpha=\epsilon_1+\epsilon_{k+1},\ldots,
\epsilon_1+\epsilon_m, \epsilon_1-\epsilon_m, \ldots,
\epsilon_1-\epsilon_{k+1}$.

\item[ii)] For the final element
$w=s_ms_{m-1}\ldots s_1$ there is only one non-degenerate
element of colength $1$ below it, namely  $ws_{\epsilon_1}$.

\item[iii)] For the final element $w=s_k\ldots s_1$ with $1\leq k \leq m-1$
there is only one
non-degenerate element of colength one below $w$, namely
 $ws_{\epsilon_1-\epsilon_{k+1}}$.
\end{enumerate}
\end{corollary}
\end{subsection}
\begin{subsection}{Final elements in $\Weyl{D}{m}$}
This section is analogous to the preceding one and we will therefore
be brief.
\begin{lemma}\label{final factorisationsD}
The products $w_k=s_k\cdots s_{m-2}s_ms_{m-1}\cdots s_1$ with
$1\leq k \leq m-2$ together with the product $w_{m-1}=s_ms_{m-1}\cdots s_1$,
the elements $w_m=s_ms_{m-2}\cdots s_1$ and
$w_{m+1}=s_{m-1}s_{m-2}\cdots s_1$
and the products $w_{m+j}=s_{m-j}s_{m-j-1}\cdots s_1$ with $j=2,\ldots, m$
are reduced expressions for the $2m$ final elements of
$\Weyl{D}{m}$. We have $w_1=w_{\emptyset}$ and $w_{2m}=1$.
\begin{proof}
Easily verified.
\end{proof}
\end{lemma}
For each integer $\ell$ with $0\leq \ell \leq 2m-2$ and $\ell \neq m-1$
there is one final element with length $\ell(w)=\ell$ while there
are two final elements of length $m-1$.
We can associate a graph to these $2m$ final elements by associating
a vertex to each final element and an edge to a pair $u,v$
if $v=s_ju$ for some $s_j$.
Conjugation by $s_m^{\prime}$ interchanges the
two final elements of length $m-1$.
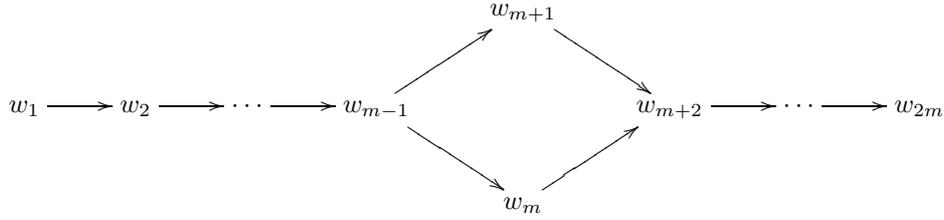
\begin{displaymath}
\begin{xy}
\xymatrix{ 
   &  &  &   & w_{m+1} \ar[dr] \\
 w_1\ar[r]  & w_2\ar[r]& \cdots   \ar[r] & w_{m-1} \ar[dr]\ar[ur] & & w_{m+2}
\ar[r]  &\cdots \ar[r] & w_{2m}\\
            &  &                         & & w_{m} \ar[ur] \\
}
\end{xy}
\end{displaymath}
We now turn to the colength $1$ elements below the final elements.
\begin{lemma}\label{colength one D}
The colength $1$ elements below the final elements are as follows:
\item[i)] There are $2m-k-1$  elements in $\Weyl{D}{m}$ of colength $1$
below the final element  $w=s_k \cdots s_{m-2} s_m \cdots s_1$
for $k\leq m-2$ and they are
$s_k \cdots \hat{s}_i \cdots s_{m-2} s_m \cdots s_1=
w s_{\epsilon_1+\epsilon_{i+1}}$  for $i=k,\ldots, m-2$,
the elements $s_1\cdots s_{m-2}\hat{s}_m s_{m-1}\cdots s_1=ws_{\epsilon_1+
\epsilon_m}$,
$s_1\cdots s_{m-2}{s}_m \hat{s}_{m-1}\cdots s_1=ws_{\epsilon_1-
\epsilon_m}$, and the elements
$s_1\cdots s_{m-2}s_ms_{m-1} \cdots \hat{s}_{m-i} \cdots s_1=
ws_{\epsilon_1-\epsilon_{m+1-i}}$ for $i=2,\ldots,m-1$.

\item[ii)] There are $m$ elements of colength $1$ below the
final element $w=s_ms_{m-1}\ldots s_1$ and they are
$ws_{\alpha}$ with $\alpha=\epsilon_1+\epsilon_m,
\epsilon_1-\epsilon_m,\ldots, \epsilon_1-\epsilon_2$.

\item[iii)] The elements in $\Weyl{D}{m}$ of colength $1$
below the final element $w=s_ms_{m-2}\cdots s_1$ are $s_{m-2}\cdots s_1=
ws_{\epsilon_1+\epsilon_m}$ and
$s_m s_{m-2} \cdots \hat{s}_i \cdots s_1=ws_{\epsilon_1-\epsilon_{i+1}}$
for $i=m-2,\ldots,1$.

\item[iv)] The elements in $\Weyl{D}{m}$ of colength $1$
below the final element $w=s_{m-1}s_{m-2}\cdots s_1$ are $s_{m-2}\cdots s_1=
ws_{\epsilon_1-\epsilon_m}$ and
$s_{m-1} s_{m-2} \cdots \hat{s}_i \cdots s_1=ws_{\epsilon_1-\epsilon_{i+1}}$
for $i=m-2,\ldots,1$.

\item[v)]  The elements in $\Weyl{D}{m}$ of colength $1$
below the final element $w=s_k\cdots s_1$ with $1 \leq k \leq m$ are
the elements $s_k \cdots \hat{s}_i \cdots s_1=ws_{\epsilon_1-\epsilon_{k+1}}$
for $i=k,\ldots,1$.
\begin{proof}
The proof is analogous to the case of $\Weyl{B}{m}$ treated in the
preceding section.
\end{proof}
\end{lemma}

Again, we now consider the elements of colength $1$ below a final element and
determine if they have degenerate or final shuffle type.
\begin{proposition}\label{colength1-typeD}
\item[i)]
The element $s_k\cdots s_{m-2}s_m\cdots\hat{s}_i\cdots s_1$ with $i<k\leq m-2$
is degenerate.

\item[ii)] The element $s_k\cdots s_{m-2}
s_m\cdots\hat{s}_i\cdots s_1$ with $i\geq k$ has an
elementary shuffle by ($s_{i+1}$) to the element $s_k\cdots s_{m-2}
s_m\cdots\hat{s}_{i+1}s_i\cdots
s_1$ if $i<m-2$ and a double shuffle (by $s_{m-1}s_m$)
to $s_k\cdots \hat{s}_{m-2}s_m\cdots s_1$
if $i=m-2$ and $k<m-2$ and an elementary shuffle (by $s_m$)
to $s_{m}\cdots s_1$ if $k=m-2=i$.

\item[iii)]
The elements $s_k\cdots s_{m-2}s_m\hat{s}_{m-1}s_{m-2}\cdots s_1$
and $s_k\cdots s_{m-2}\hat{s}_ms_{m-1}\cdots s_1$
are degenerate.
\item[iv)] For $m-2\geq i>k$ the element
$s_k\cdots\hat{s}_i\cdots s_{m-2} s_m\cdots s_1$ has an elementary shuffle
(by $s_i$)  to
the element $s_k\cdots \hat{s}_{i-1}\cdots s_{m-2} s_m\cdots s_1$.

\item[v)] The element $s_m\cdots \hat{s}_i\cdots s_1$
is degenerate if $m-1\leq i \leq 1$.

\item[vi)] The element $s_m s_{m-2} \cdots \hat{s}_i \cdots s_1$ with $m-2\geq i
>1$ is degenerate.
\begin{proof}
The proof is analogous to the $\Weyl{B}{m}$ case and is omitted.
\end{proof}
\end{proposition}
\begin{corollary}\label{colength1-typeD-cor}
The non-degenerate elements of colength 1 below the final ones are as follows.
\item[i)] For a final element $w=s_k\ldots s_{m-2}s_m\ldots s_1$
with $k\leq m-2$
the only non-degenerate elements of colength $1$ below $w$ are
$ws_{\alpha}$ with $\alpha=\epsilon_1+\epsilon_{k+1},\ldots,
\epsilon_1+\epsilon_{m-1}, \epsilon_1-\epsilon_{m-1}, \ldots,
\epsilon_1-\epsilon_{k+1}$.
\item[ii)] For the final element $s_ms_{m-1}\ldots s_1$
there are two non-degenerate elements of colength $1$ below it,
namely $ws_{\epsilon_1+\epsilon_m}$ and $ws_{\epsilon_1-\epsilon_m}$

\item[iii)] For the final element $w=s_ms_{m-2}\ldots s_1$ there is only one
non-degenerate element of colength~$1$ below it,
namely $ws_{\epsilon_1+\epsilon_m}$.

\item[iv)] For the final element $w=s_{m-1}s_{m-2}\ldots s_1$ there is only
one non-degenerate element of colength~$1$ below it, namely
$ws_{\epsilon_1-\epsilon_m}$.

\item[v)] For the final element $w=s_k \ldots s_1$ with $1\leq k \leq m-2$
there is only one non-degenerate element of colength~$1$ below it, namely
$ws_{\epsilon_1-\epsilon_{k+1}}$.
\end{corollary}

\end{subsection}
\begin{subsection}{Twisted final elements in $\Weylp{D}{m}$}
The twisted final elements are of the form $ws_m'$ with $w$ as in
Lemma \ref{final factorisationsD}. Similarly, the elements of colength $1$
below a twisted final element $ws_m'$ are of the form $u s_m'$ with
$u$ a colength $1$ element below $w$ as described in Lemma
\ref{colength one D}. We have to analyze whether these elements
are degenerate or have final shuffle type. We omit the analogue
of Proposition \ref{colength1-typeD} and formulate immediately the analogue of
Corollary \ref{colength1-typeD-cor}
\begin{corollary}
The non-degenerate elements of colength 1 below the final ones are as follows.
\item[i)] For a twisted final element $ws_m'$ with
$w=s_k\ldots s_{m-2}s_m\ldots s_1$
and $k\leq m-2$
the only non-degenerate elements of colength $1$ below $ws_m'$ are of the
form $us_m'$ with $u$ equal to
$ws_{\alpha}$ with $\alpha=\epsilon_1+\epsilon_{k+1},\ldots,
\epsilon_1+\epsilon_{m}$, and $\epsilon_1-\epsilon_{m}, \ldots,
\epsilon_1-\epsilon_{k+1}$.
\item[ii)] For the twisted final element $ws_m'$ with $w=s_ms_{m-1}\ldots s_1$
there are two non-degenerate elements of colength $1$ below it,
namely corresponding to
$ws_{\epsilon_1+\epsilon_m}$ and $ws_{\epsilon_1-\epsilon_m}$

\item[iii)] For the twisted final element $ws_m'$ with
$w=s_ms_{m-2}\ldots s_1$ there is only one
non-degenerate element of colength~$1$ below it,
namely corresponding to $ws_{\epsilon_1+\epsilon_m}$.

\item[iv)] For the final element $ws_m'$ with
$w=s_{m-1}s_{m-2}\ldots s_1$ there is only
one non-degenerate element of colength~$1$ below it, namely
corresponding to $ws_{\epsilon_1-\epsilon_m}$.

\item[v)] For the final element $ws_m'$ with
$w=s_k \ldots s_1$ with $1\leq k \leq m-2$
there is only one non-degenerate element of colength~$1$ below it, namely
corresponding to $ws_{\epsilon_1-\epsilon_{k+1}}$.

\end{corollary}

\end{subsection}
\end{section}

\begin{section}{Pieri's formula and the cycle classes of the strata}

  Since the strata in our case, unlike the case of abelian varieties, 
are (almost) linearly
  ordered we can fruitfully apply a Pieri type formula to get a formula for the
  classes of $\Vc_\nu$. The appropriate formula is the Pieri formula of Pittie
  and Ram (\cite{pittie99::ae+pieri+cheval+k+g+b}).  There is a small problem in
  that the result only applies when we start with a $G$-torsor over a connected
  semi-simple group $G$ and in our case the structure group is the disconnected
  group $\Ogrp(n)$. The resolution of this problem differs somewhat in the two
  cases of even or odd $n$ so part of the discussion is postponed to the
  separate discussions for the two cases. In any case, the Pittie-Ram formula
  expresses the intersection product of the cycle class of a stratum with a
  first Chern class in terms of cycle classes of strata of one dimension less.
  We have to use the formula on the flag space and then project it down. The
  precise details of the Pieri formula differ enough between the odd and even
  dimensional cases to make separate discussions in the two cases.
Throughout we shall assume the versality assumption made in Section \ref{localstructure}:
we assume that we have a family $f: X \to S$ of $N$-marked K3 surfaces 
(where S may be an algebraic stack) such that
$S$ be smooth over ${\FF}_p$ at  all $s$ and that the composed map 
$T_sS\to {\rm Hom}(H^0(X_s,\Omega^2_{X_s}),P)$ be  surjective.

\begin{subsection}{The odd-dimensional case}

We now assume $n=2m+1$. Now, ${\Ogrp}(2m+1)={\SO}(2m+1)\times\{\pm1\}$ 
and hence an ${\Ogrp}(2m+1)$-torsor is the same thing as one ${\SO}(2m+1)$-torsor and one double
cover. Thus the problem mentioned above is resolved by considering instead the
${\SO}(2m+1)$-torsor. In concrete terms this means replacing our $F$-zip vector
bundle $H$ by $H\Tensor\det(H)$.

We want to apply the Pieri formula to the two complete flags that we have on the
flag space $\cF_n$.  If we let $\lambda = \sum_in_i\ell_i$, where $\ell_i =
c_1(E_i/E_{i-1})$ for $1\leq i \leq m$ is the first Chern class corresponding
to the root $\epsilon_i$.  The starting point for the Pieri formula is the
following construction: Given a sequence $z=(z_1,\dots,z_m)$ of cohomology
classes (of fixed degree) and a weight vector
$\lambda=\sum_{i=1}^mn_i\epsilon_i$ (in the weight lattice of type $B_m$) we
define $z^\lambda:=\sum_in_iz_i$. We shall apply this to
$x=(\ell_1,\dots,\ell_m)$ and $y=(k_1,\dots,k_m)$, where
$k_i:=c_1(G_i/G_{i-1})$ and then, for suitable $\lambda$, we shall consider
$x^\lambda$ and $y^{w\lambda}$. However, the elements of $x$ and of $y$ span the
same subgroup of the cohomology so we can also write $y^{w\lambda}$ as
$x^{\lambda '}$ for a suitable $\lambda '$. Clearly the association $\lambda
\mapsto \lambda '$ is a linear operator on the weight lattice. It is easily seen
that just as for the symplectic case it is given by $\lambda '= pw_\emptyset
w(\lambda)$. From this point on we shall only be considering elements of the
form $x^\mu$ and for simplicity we shall write them just as $\mu$. The Pieri
formula (see the proof of \cite[Thm 10.1]{ekedahl10::cycle+class+e+o+strat} for
details) now takes the form
\begin{displaymath}
(1-pw_\emptyset w)(\lambda)[\Uc_w]=
-\sum_{\ell(ws_\alpha)=
\ell(w)-1}\langle\alpha^{\vee},\lambda\rangle[\Uc_{ws_\alpha}].
\end{displaymath}
The term $1-pw_\emptyset w$ is viewed as an element of the group ring
${\QQ}[\Weyl{B}{m}]$ acting on the roots $\ell_i$ and the sum is over roots
$\alpha$ such that the length $\ell(ws_{\alpha})$ is one less than the length of
$w$. Moreover, $\alpha^{\vee}$ is the usual coroot defined by $\alpha$.
To obtain a formula for the multiplication of $[\Uc_w]$ by a given line bundle
$\rho$ we have to solve the equation $(1-pw_\emptyset w)(\lambda)=\rho$. If we
put $v:=w_\emptyset w$ and if we let $c$ be the smallest positive integer such
that $v^c(\rho)=s\rho$ for some $s \in \{\pm1\}$ then a solution is given by
\begin{displaymath}
\lambda = \frac{1}{1-sp^c}\sum_{i=0}^{c-1}p^iv^i(\rho).
\end{displaymath}
We carry this out with $\rho=\ell_1=\lambda_1$, the first Chern class of the
Hodge bundle, such that we obtain a formula for $\lambda_1 [\Uc_w]$. We shall
call $c$ the {\sl reduced orbit length} and say that the orbit is
{\sl even} or {\sl odd} according to as $s$ is $+1$ or $-1$.  Then
we push down to the moduli space. The degenerate strata push down to zero and
the non-degenerate to a power of $p$ times the push down of a final stratum.

The final elements in this case are of the form $w_k=s_k\cdots s_m\cdots s_1$
with $1 \leq k<m$ and $w_{k+m}=s_{m-k}s_{m-k-1} \cdots s_1$ for $0 \leq k \leq
m-1$ and $w_{2m}=1$. Note that $w_{\emptyset}= w_1$ with this usage.  The
corresponding final strata on the flag space are $\Uc_{w_k}$ with
$k=1,\ldots,2m$ with corresponding strata $\Vc_{w_k}$ on the moduli space.
For a final element we denote the canonical map
$\Uc_w \to \Vc_w$ by $\pi_w$ and its degree by $\deg(\pi_w)$.
\begin{remark}
The strata $\Vc_{w_k}$ for $1\leq k \leq m$ are the strata corresponding to
finite height equal to $\geq k$, the stratum $\Vc_{w_{m+1}}$ is the supersingular
stratum and the stratum $\Vc_{w_{k+m}}$ corresponds to Artin invariant 
$\leq m+1-k$ for $1\leq k \leq m$.
\end{remark}
\begin{theorem}\label{Odd cycle classes}
  The cycle classes of the final strata $\Vc_w$ on the base $S$ are powers
  of $\lambda_1$ times polynomials in $p$ 
  given by
\begin{eqnarray*}
{\rm i)} \quad
[\Vc_{w_k}] &=& (p-1)(p^2-1)\cdots(p^{k-1}-1) \lambda_1^{k-1} \quad
\hbox{\rm if $1\leq k\leq m$,}\\
{\rm ii)} \quad [\Vc_{w_{m+1}}] &=&\frac{1}{2} (p-1)(p^2-1)\cdots(p^{m}-1)
\lambda_1^{m},\\
{\rm iii)} \quad
[\Vc_{w_{m+k}}] &=&\frac{1}{2}
\frac{(p^{2k}-1)(p^{2(k+1)}-1)\cdots(p^{2m}-1)}{(p+1)\cdots(p^{m-k+1}+1)}
\lambda_1^{m+k-1} \quad
\hbox{\rm if $2\leq k\leq m$.}
\end{eqnarray*}
\begin{proof}
We start with a final element $w$ of the form $w_k=s_k\cdots s_m\cdots s_1$
with $1 \leq k<m$.  The colength $1$ elements $w_ks_{\alpha}$
that are not degenerate correspond to the $2m-2k$ elements
$\alpha_1= \epsilon_{1}+\epsilon_{k+1},
\ldots, \alpha_{m-k}= \epsilon_{1}+\epsilon_{m}, \alpha_{m-k+1}=
\epsilon_1-\epsilon_m, \ldots, \alpha_{2m-2k}=\epsilon_1-\epsilon_{k+1}$.
These are the only elements that will contribute to the push down.
Note that we have $w_ks_{\alpha_1}=w_{k+1}$, again a final element.
For the element $v=w_{\emptyset}w$ we have $v^j(1)=2m+1-k+j$ for
$j=1,\ldots,k-1$ which means that the reduced orbit length is $k$ and the orbit
even. We thus find that
$$
(1-pv) \lambda_1 = \sum_{i=0}^{k-1} p^i v^i (\ell_1)=
\ell_1 - \sum_{i=1}^{k-1} p^i \ell_{k+1-i}.
$$
Therefore the Pieri formula gives
\begin{eqnarray*}
(p^{k}-1) \lambda_1 [\Uc_{w_k}] \equiv
&\sum_{j=1}^{m-k}  (\epsilon_1+\epsilon_{k+j},
\ell_1-\sum_{i=1}^{k-1}p^i \ell_{k+1-i})
[\Uc_{ws_{\alpha_j}}] +\\
&\sum_{j=1}^{m-k}  (\epsilon_1-\epsilon_{m+1-j},
\ell_1-\sum_{i=1}^{k-1}p^i \ell_{k+1-i})
[\Uc_{ws_{\alpha_{m-k+j}}}],
\end{eqnarray*}
where $\equiv$ means that we count modulo degenerate strata.
Pushing it down annihilates the classes of the degenerate strata because these
loose dimension and yields
\begin{eqnarray*}
(p^k-1) \lambda_1 [\Vc_{w_k}] \deg (\pi_{w_k}) & = &
\sum_{j=1}^{2m-2k} [\Vc_{w_{k+1}}]
\deg(\pi_{ws_{\alpha_j}}) \\
&=& (1+p+\ldots +p^{2m-2k-1}) [\Vc_{ws_{k+1}}] \deg (\pi_{w_{k+1}})
\\
\end{eqnarray*}
since the $w_ks_{\alpha_j}$ for $j=2,\ldots, 2m-2k$ are shuffles of
$w_{k+1}=ws_{\alpha_1}$ which map to $\Vc_{w_{k+1}}$ with degree $p^{j-1}$.
By Lemma \ref{DegreelemmaB} we have
$\deg(\pi_{w_k})/\deg(\pi_{w_{k+1}})=p^{2m-2k-1}+\ldots+1$ and get
$[\Vc_{w_{k+1}}]=(p^k-1) [\Vc_{w_k}]$ for $k=1,\ldots, m-1$. Since
$[\Vc_{w_1}]=1$ part i) follows.

For part ii) we note that there is only one
non-degenerate element of colength $1$, namely $ws_{\alpha}=w_{m+1}$ and
it corresponds to $\alpha = \epsilon_1$ with $\alpha^{\vee}=2\epsilon_1$.
Note that $v=[m+2,2m+1,2,3,\ldots,m-1]$ and
$\sum_{i=0}^{m-1} p^iv^i(\ell_1)=\ell_1-\sum_{i=1}^{m-1}p^i\ell_{m+1-i}$.
This gives
$$
(p^m-1) [\Vc_{w_m}] \deg (\pi_{w_m})= 2 \,  [\Vc_{w_{m+1}}] \deg (\pi_{w_{m+1}})
$$
and we observe that $\deg (\pi_{w_m})=1=\deg (\pi_{w_{m+1}})$.
This proves ii).
For the case iii) we consider a final element
$w_{k+m}=s_{m-k}s_{m-k-1}\cdots s_1$ with $k\geq 1$.
There is only one non-degenerate element
$w_{k+m}s_{\alpha}$ of colength $1$, namely $w_{k+1+m}$
with $\alpha=\epsilon_{1}-\epsilon_{m+1-k}$.

The element $v=w_{\emptyset}w_{m+k}=[m,2m+1,2,3,\ldots]$ has an odd orbit of
reduced orbit length $m+1-k$ and thus $v^{m+1-k}\lambda_1=-\lambda_1$ so that
$\sum_{i=0}^{m-k} p^iv^i \ell_1 = \ell_1+\sum_{i=1}^{m-k} p^i \ell_{m+2-k-i}$
and the Pieri formula gives
$$
(p^{m-k}+1)\lambda_1 [\Vc_{m+k}]\deg (\pi_{w_{m+k}})=
(p-1) [\Vc_{w_{m+k+1}}] \deg (\pi_{w_{m+k+1}}).
$$
Here we have
$\deg (\pi_{w_{m+k+1}})/\deg (\pi_{w_{m+k}})= p^{2k-1}+\cdots+1$
which gives
$(p^{m-k}+1) \lambda_1[\Vc_{m+k}]= (p^{2k}-1) [\Vc_{w_{m+k+1}}]$.
This proves the formulas.

That the formulas are up to a factor $1/2$ polynomials in $\ZZ[\lambda_1,p]$
is clear for cases i) and ii) and follows from the next remark for case iii).
\end{proof}
\end{theorem}
\begin{remark}
The formula for case iii)  can also be written as
\begin{displaymath}
[\Vc_{w_{m+k}}] =\frac{1}{2} \left( \prod_{j=1}^{m+1-k} (p^j-1) \right)
\genfrac{[}{]}{}{}{m}{m+1-k}_{p^2} \lambda_1^{m+k-1},
\end{displaymath}
where $\genfrac{[}{]}{}{}{n}{i}_q$ is the usual $q$-binomial coefficient.
\end{remark}

\begin{remark}
Theorem \ref{examplethm} in the Introduction is the special case 
where we take $S$
equal to the moduli space of of polarized K3 surfaces of degree 
$d$, prime to $p$ and where $m=10$. The versality condition is verified
in a standard way (and essentially the same as the proof for the case of
elliptic K3 surfaces with a section below).
\end{remark}

\end{subsection}
\begin{subsection}{The even-dimensional case}
The reduction to an ${\SO}$-torsor in the even case is more involved than in the
odd case. To begin with if we have an ${\Ogrp}(2m)$-torsor $P \to X$ we get a
double cover $Y:=P/{\SO}(2m) \to X$ and the quotient map $P \to Y$ is an
${\SO}(2m)$-torsor. However, in order to have a Bruhat cell decomposition of $P
\to Y$ (which is necessary even to formulate the Pieri formula) we need a
reduction of the structure group to $B$ (a Borel subgroup of ${\SO}(2m)$). This we
get from our original setup in the following way: We can find a subgroup
$B'\subset {\Ogrp}(2m)$ containing $B$ as a subgroup of index $2$ and we assume
that we have a $B'$-torsor $Q \to X$ which then gives rise to a $B$-torsor $Q
\to Q/B=Y$. If we look at the corresponding $G/B$-fibrations 
(where $G={\SO}(2m)$)
we get a commutative diagram
\begin{displaymath}
\begin{CD}
Q\times_BG/B @>>> Q\times_{B'}G/B\\
@VVV    @VVV\\
Y @>>>  X
\end{CD}
\end{displaymath}
and we get a Pieri formula to $Q\times_BG/B \to Y$ and then push it down to $
Q\times_{B'}G/B$. What happens during this pushdown is the following: The class
$\lambda_1$ (which is the only class for which we shall use the Pieri formula)
is the pullback of a class on $Q\times_{B'}G/B$ so by the projection formula it
can be moved out of the push down. We get a ``Bruhat decomposition'' also of
$Q\times_{B'}G/B$ but the strata now corresponds to $B'$-orbits of $G/B$. Such
an orbit is a union of one or two $B$-orbits depending on whether or not an
element in $B' \setminus B$ fixes the $B$-orbit or not. Hence, the projection
$Q\times_BG/B \to Q\times_{B'}G/B$ maps Bruhat strata to Bruhat strata and two
strata $\Uc_w$ and $\Uc_{w'}$ in $Q\times_BG/B$ are mapped to the same stratum
precisely when either $w'=w$ or $w'=s'_mws'_m$. In our specific case the
$B'$-bundle arises by starting with a $G'$-torsor ($G'={\Ogrp}(2m)$) $P \to X$ and
then pulling it back along $P\times_{G'}G'/B'$ (note that $G'/B'=G/B$) where
this pullback has a canonical reduction of its structure group to $B'$.

In flag terms we have the following description. Let $E$ be a quadratic vector
bundle of rank $2m$ over $X$. The pairing on it induces an isomorphism
$\det(E)\Tensor\det(E) \riso \sO_X$ and hence gives a double cover $Y \to X$. We
can also consider the {\sl almost complete} flag space $\sF \to X$ of
self dual flags $0\subset E_1 \subset E_2 \subset \cdots \subset E_{m-1}\subset
E_{m+1}\subset\cdots\subset E_{2m}=E$ with $\dim E_i=i$. The fibre product
$\sF':=Y\times_X\sF$ has the explicit description as the space of complete self
dual flags $0\subset E_1 \subset E_2 \subset \cdots \subset E_{m-1}\subset E_m
\subset E_{m+1}\subset\cdots\subset E_{2m}=E$ and the double cover involution on
$Y$ induces the operation on such flags which replaces $E_m$ by the other
totally isotropic $m$-dimensional subspace $E_{m-1}\subset E'_m \subset
E_{m+1}$. The fibre product $\sF'':=\sF'\times_X\sF'$ consisting of pairs
$(E_\bullet,F_\bullet)$ of complete self dual flags split up in
two components: one is $\sF''_0$, where $\dim(E_m\cap F_m)\equiv m \bmod 2$,
and the other one is $\sF''_1$, where $\dim(E_m\cap F_m)\equiv m+1 \bmod 2$. The group
$\ZZ/2\times\ZZ/2$ acts on $\sF'\times_X\sF'$, a group factor acting on the
corresponding factor of $\sF'\times_X\sF'$. The elements $(1,0)$ and $(0,1)$
permutes the two components and $(1,1)$ preserves them. Each element $w \in
\Weylp{D}{m}$ gives a stratum of $\sF''$ consisting of the flags in relative
position $w$. When $w \in \Weyl{D}{m}$, the stratum lies in $\sF''_0$ and when
$w \in \Weyl{D}{m}s'_m$ it lies in $\sF''_1$. The group element $(1,0)$ then
takes the stratum of $w$ to that of $ws'_m$ and the element $(0,1)$ takes $w$ to
that of $s'_mw$.
\end{subsection}
\begin{subsubsection}{The untwisted even case}

The first step in getting to a Pieri formula is to identify the linear map that
takes $\lambda$ to $\lambda '$. This time it is \emph{not} given by $pw_\emptyset w$
as $w_\emptyset(m)=m+1$ which does not have the desired effect. Instead we have
to use the linear map $pw'_\emptyset w$ because $w'_\emptyset(m)=m$.  This means
that Pieri's formula takes the form
\begin{displaymath}
(1-pw'_\emptyset w)(\lambda)[\Uc_w]=
-\sum_{\ell(ws_\alpha)=
\ell(w)-1}\langle\alpha^{\vee},\lambda\rangle[\Uc_{ws_\alpha}]
\end{displaymath}
Recall that the final elements are the $2m$ elements $w_k=s_k\ldots
s_{m-2}s_m\ldots s_1$ for $k=1,\ldots, m-2$, $w_{m-1}=s_ms_{m-1}\ldots s_1$,
$w_m=s_ms_{m-2}\ldots s_1$ and $w_{m+j}=s_{m-j}\ldots s_1$ for $j= 1,\ldots,m-1$
and $w_{2m}=1$. Moreover, there is an automorphism of $\Weyl{D}{m}$
interchanging $w_m$ and $w_{m+1}$ given by conjugation by $s'_m$.
\begin{theorem}\label{Untwisted classes}
The cycle classes of the final strata $\Vc_w$ for final
$w_j \in \Weyl{D}{m}$ on the base $S$ are
powers of $\lambda_1$ times polynomials
in $p$ given by
\begin{eqnarray*}
{\rm i)} \quad
[\Vc_{w_k}] &=& (p-1)(p^2-1)\cdots(p^{k-1}-1) \lambda_1^{k-1} \quad
\hbox{\rm if $k\leq m-1$, } \\
{\rm ii)} \quad [\Vc_{w_{m+1}}]
&=& (p-1)(p^2-1)\cdots(p^{m-1}-1)
\lambda_1^{m-1},\\
{\rm iii)} \quad
[\Vc_{w_{m+k}}] &=&
\frac{1}{2}
\frac{\prod_{i=1}^{m-1}(p^i-1)\prod_{i=m-k+2}^m(p^i+1)}
{\prod_{i=1}^{k-2}(p^i+1)\prod_{i=1}^{k-1}(p^i-1)}
\lambda_1^{m+k-2} \quad
\hbox{\rm if $2\leq k\leq m$. }
\end{eqnarray*}
Furthermore, we have that $\Vc_{w_m}=\emptyset$.
\begin{proof}
  Let $w_k=s_k\ldots s_{m-2}s_m\ldots s_1$ be a final element with $1 \leq k
  \leq m-2$. There are $2m-2k-2$ non-degenerate elements of colength $1$ under
  $w_k$ and they are of the form $ws_{\alpha}$ with
  $\alpha_{1}=\epsilon_1+\epsilon_{k+1},\ldots,\alpha_{m-k-1}=\epsilon_1+\epsilon_{m-1}$
  and $\alpha_{m-k}=\epsilon_1-\epsilon_{m-1},\ldots,
  \alpha_{2m-2k-2}=\epsilon_1-\epsilon_{k+1}$.  We find that
  $v:=w'_{\emptyset}w_k$ has an even orbit of reduced orbit length $k$ and
$$
(1-pv)\lambda_1=\ell_1- \sum_{i=1}^{k-1} p^i \ell_{k+1-i}.
$$
Therefore the Pieri formula gives
\begin{eqnarray*}
(p^{k}-1) \lambda_1 [\Uc_{w_k}]  \equiv
&\sum_{j=1}^{m-k-1}  (\epsilon_1+\epsilon_{k+j},
\ell_1-\sum_{i=1}^{k-1}p^i \ell_{k+1-i})
[\Uc_{ws_{\alpha_j}}] +\\
&\sum_{j=m-k}^{2m-2k-2}  (\epsilon_1-\epsilon_{2m-k-1-j},
\ell_1-\sum_{i=1}^{k-1}p^i \ell_{k+1-i})
[\Uc_{ws_{\alpha_{j}}}],
\end{eqnarray*}
where $\equiv$ means again that we work modulo degenerate strata.
Pushing it down annihilates the classes of the degenerate strata
and yields
\begin{eqnarray*}
(p^k-1) \lambda_1 [\Vc_{w_k}] \deg (\pi_{w_k}) & = &
\sum_{j=1}^{2m-2k-2} [\Vc_{w_{k+1}}]
\deg(\pi_{ws_{\alpha_j}}) \\
&=& (1+\ldots +p^{m-k-2}+p^{m-k} +\ldots +p^{2m-2k-2})
[\Vc_{w_{k+1}}] \deg (\pi_{w_{k+1}}), \\
\end{eqnarray*}
since the $w_ks_{\alpha_j}$ for $j=1,\dots,m-k-1$ are shuffles of
$w_{k+1}=ws_{\alpha_{0}}$ for which $\Uc_{w_ks_{\alpha_j}}$ maps to
$\Uc_{w_{k+1}}$ with degree $p^{j-1}$, while for $j=m-k,\ldots, 2m-2k-2$ we get
degree $p^{j}$.  By Lemma \ref{DegreelemmaB} we have
$\deg(\pi_{w_k})/\deg(\pi_{w_{k+1}})=p^{2m-2k-2}+\ldots+p^{m-k}+p^{m-k-2}+\dots+1$
and hence get $[\Vc_{w_{k+1}}]=(p^k-1)\lambda_1[\Vc_{w_k}]$ for $k=1,\ldots,
m-1$. Since $[\Vc_{w_1}]=1$ part i) follows.

For part ii) we consider the final element $w=s_ms_{m-1}\ldots s_1$ and see that
$v:=w'_{\emptyset} w$ has an even orbit of reduced orbit length $m-1$ and that
$(1-pv)\lambda_1=\ell_1-\sum_{i=1}^{m-2} p^i \ell_{m-i}$.  In this case there
are two non-degenerate elements of colength $1$ namely $ws_{\alpha}$ with
$\alpha$ being equal to $\epsilon_1+\epsilon_m$ and $\epsilon_1-\epsilon_m$
respectively. Applying the Pieri formula and pushing down first to the
unoriented flag space and then to the moduli space we get
$$
(p^{m-1}-1)\lambda_1[\Vc_{m-1}] \deg(\pi_{w_{m-1}})=
[\Vc_{w_{m}+w_{m+1}}] \deg(\pi_{w_m}).
$$
By Lemma \ref{DegreelemmaB} $\deg(\pi_{w_{m-1}})=\deg(\pi_{w_m})=1$ which gives
ii) after pushing down.

For part iii) we consider the element $w_m=s_ms_{m-2}\ldots s_1$.  The element
$v:=w'_{\emptyset} w$ has an odd orbit of reduced orbit length $m$ and we have
$$
(1-pv)\lambda_1= \ell_1-p\ell_m+\sum_{i=2}^{m-1} p^i \ell_{m+1-i}.
$$
There is now only one non-degenerate element of colength $1$,
namely $w_ms_{\alpha}$ with $\alpha=\epsilon_1+\epsilon_m$.
We get
$$
(p^m+1)\lambda_1 [\Vc_{w_m+w_{m-1}}]\deg(\pi_{w_{m}})= (p-1) [\Vc_{w_{m+2}}]
\deg(\pi_{w_{m+2}}).
$$
Again by Lemma \ref{DegreelemmaD} we have
$\deg(\pi_{w_{m}})=\deg(\pi_{w_{m+2}})=1$.
Now take $w=w_{m+j}=s_{m-j}\ldots s_1$ with $j\geq 2$. The element
$v:=w'_{\emptyset}w$ has an odd orbit of reduced orbit length $m-j+1$. We get
$$
(1-pv) \lambda_1= \ell_1+\sum_{i=1}^{m-j} p^i \epsilon_{m+2-j-i}.
$$
There is again only one non-degenerate element $ws_{\alpha}$ with
$\alpha=\epsilon_1-\epsilon_{m+1-j}$. Therefore
$\langle \alpha^{\vee},\lambda\rangle= (1-p^{j-1})$. We find
$$
(p^{m+1-j}+1) \lambda_1 [\Vc_{w_{m+j}}]
\deg(\pi_{w_{m+j}})=(p^{j-1}-1) [\Vc_{w_{m+j}}] \deg(\pi_{w_{m+j+1}}).
$$
By Lemma \ref{DegreelemmaD} we have
$\deg(\pi_{w_{m+j+1}})\deg(\pi_{w_{m+j}})=
p^{2j-2}+\cdots+2p^{j-1}+\cdots+1$. Using the
factorization
$p^{2j-2}+\cdots+2p^{j-1}+\cdots+1=(p^{j-1}+1)(p^{j-1}+p^{j-2}+\cdots+1)=
(1+p^{j-1})/(1-p^j)$ and
iterating we get the formula. As there is only a small number cases ($m\leq 10$), the fact
that one gets polynomials is most easily verified by explicit computation (one
could also use \cite[Prop.~13.2]{ekedahl10::cycle+class+e+o+strat}).
\end{proof}
\end{theorem}
\begin{remark}
Similarly to the odd case we can write the third formula as
\begin{displaymath}
\frac{p^{k-1}+1}{p^m-1}\left( \prod_{i=1}^{m+1-k}(p^i-1) \right)
\left[\begin{array}{c}
m\\
k-1
\end{array}\right]_{p^2},
\end{displaymath}
though this does not make it visibly a polynomial in $p$.
\end{remark}
\end{subsubsection}
\begin{subsubsection}{The twisted even case}

We now turn to the twisted even-dimensional case. Going back to the previous
notation we have the space $\sF''$ of pairs of complete flags of $H$ where $H$
now is the de Rham cohomology of the universal K3 surface over our moduli space. We have a disjoint
decomposition $\sF''=\sF''_0\sqcup \sF''_1$, where $\sF''_0$ is the
$G/B$-fibration over $\sF'$ with structure group $B$ and where we consequently
have a Pieri formula. However, the $F$-zip structure on $H$ gives a section of the
projection (on the first factor) $\sF'' \to \sF'$ which is contained completely
in $\sF''_1$. In order to get a section along which we can pullback a Pieri
formula on $\sF''_0$ we must compose with the isomorphism $\sF''_1 \riso
\sF''_0$ obtained by applying the involution of $\sF'$ acting on the second
factor (say) of $\sF''$. This extra involution implies that the linear map
$\lambda \mapsto \lambda '$ is now given by $pw_\emptyset w$ (and \emph{not} by
$pw'_\emptyset w$ as in the untwisted case). Apart from that the argument of the
Pieri formula proceeds along lines very similar to the untwisted case.

Recall that the final elements are the $2m$ elements $w_k=s_k\ldots
s_{m-2}s_m\ldots s_1$ for $k=1,\ldots, m-2$, $w_{m-1}=s_ms_{m-1}\ldots s_1$,
$w_m=s_ms_{m-2}\ldots s_1$ and $w_{m+j}=s_{m-j}\ldots s_1$ for $j= 1,\ldots,m-1$
and $w_{2m}=1$.
\begin{theorem}\label{Twisted classes}
The cycle classes of the final strata $\Vc_w$ for twisted final elements
$w_j \in \Weyl{D}{m} s_m^{\prime}$ on the base $S$ are
powers in $\lambda_1$ with coefficients 
that are polynomials in $p$ given by
\begin{eqnarray*}
{\rm i)} \quad
[\Vc_{w_k}] &=& (p-1)(p^2-1)\cdots(p^{k-1}-1) \lambda_1^{k-1} \quad
\hbox{\rm if $k\leq m-1$, } \\
{\rm ii)} \quad [\Vc_{w_m}]
&=& (p-1)(p^2-1)\cdots(p^{m}-1)
\lambda_1^{m-1},\\
{\rm iii)} \quad
[\Vc_{w_{m+k}}] &=&
\frac{1}{2}
\frac{\prod_{i=1}^{m}(p^i-1)\prod_{i=m-k+2}^{m-1}(p^i+1)}
{\prod_{i=1}^{k-1}(p^i+1)\prod_{i=1}^{k-2}(p^i-1)}
\lambda_1^{m+k-2} \quad
\hbox{\rm if $2\leq k\leq m$.}
\end{eqnarray*}
Furthermore, we have $\Vc_{w_{m+1}}=\emptyset$.
\begin{proof}
Let $w_k=s_k\ldots s_{m-2}s_m\ldots s_1$ be a final element with $1 \leq k \leq
m-2$. There are $2m-2k$ non-degenerate elements of colength $1$ under $w_k$ and
they are of the form $ws_{\alpha}$ with
$\alpha_{1}=\epsilon_1+\epsilon_{k+1},\ldots,\alpha_{m-k}=\epsilon_1+\epsilon_{m}$
and $\alpha_{m-k+1}=\epsilon_1-\epsilon_{m},\ldots,
\alpha_{2m-2k}=\epsilon_1-\epsilon_{k+1}$.  We find that $v:=w_{\emptyset}w_k$
has an even orbit of reduced orbit length $k$ and
$$
(1-pv)\lambda_1=\ell_1- \sum_{i=1}^{k-1} p^i \ell_{k+1-i}
$$
Therefore the Pieri formula gives
\begin{eqnarray*}
(p^{k}-1) \lambda_1 [\Uc_{w_k}]  \equiv
&\sum_{j=0}^{m-k}  (\epsilon_1+\epsilon_{k+j},
\ell_1-\sum_{i=1}^{k-1}p^i \ell_{k+1-i})
[\Uc_{ws_{\alpha_j}}] +\\
&\sum_{j=m-k+1}^{2m-2k}  (\epsilon_1-\epsilon_{2m-k+1-j},
\ell_1-\sum_{i=1}^{k-1}p^i \ell_{k+1-i})
[\Uc_{ws_{\alpha_{-j}}}],
\end{eqnarray*}
where $\equiv$ means again that we work modulo degenerate strata.
Pushing it down annihilates the classes of the degenerate strata
and yields
\begin{eqnarray*}
(p^k-1) \lambda_1 [\Vc_{w_k}] \deg (\pi_{w_k}) & = &
\sum_{j=1}^{2m-2k} [\Vc_{w_{k+1}}]
\deg(\pi_{ws_{\alpha_j}}) \\
&=& (1+\ldots +2p^{m-k-1} +\ldots +p^{2m-2k-2})
[\Vc_{w_{k+1}}] \deg (\pi_{w_{k+1}}), \\
\end{eqnarray*}
since the $w_ks_{\alpha_j}$ for $j=1,\ldots, m-k$ are shuffles of
$w_{k+1}=ws_{\alpha_{k+1}}$ for which $\Uc_{w_ks_{\alpha_j}}$ maps to
$\Uc_{w_{k+1}}$ with degree $p^{j-1}$, while for $j=m-k+1,\ldots, 2m-2k$ we get
degree $p^{j-2}$. By Lemma \ref{DegreelemmaD} we have
$\deg(\pi_{w_k})/\deg(\pi_{w_{k+1}})=p^{2m-2k-2}+\ldots +2p^{m-k-1}+\dots+1$ and
hence get $[\Vc_{w_{k+1}}]=(p^k-1) [\Vc_{w_k}]$ for $k=1,\ldots, m-1$. Since
$[\Vc_{w_1}]=1$ part i) follows.
For part ii) we consider the final element $w=s_ms_{m-1}\ldots s_1$ and see that
$v:=w_{\emptyset} w$ has an even orbit of reduced orbit length $m-1$ and that
$(1-pv)\lambda_1=\ell_1-\sum_{i=1}^{m-2} p^i \ell_{m-i}$.  In this case there
are two non-degenerate elements of colength $1$ namely $ws_{\alpha}$ with
$\alpha_1=\epsilon_1+\epsilon_m$ and $ \alpha_2=\epsilon_1-\epsilon_m$. Applying
the formula we get
$$
(p^{m-1}-1)\lambda_1[\Vc_{m-1}] \deg(\pi_{w_{m-1}})=
[\Vc_{ws_{\alpha_1}}] \deg(\pi_{w_m})+[\Vc_{ws_{\alpha_2}}] \deg({w_{m+1}}).
$$
By Lemma \ref{DegreelemmaD} $\deg(\pi_{w_{m-1}})=\deg(\pi_{w_m})=1$ which gives
ii) after pushing down.

For part iii) we consider the element $w_m=s_ms_{m-2}\ldots s_1$.  The element
$v:=w_{\emptyset} w$ has an even orbit of reduced orbit length $m$ and we have
$$
(1-pv)\lambda_1= \ell_1+p\ell_m-\sum_{i=2}^{m-1} p^i \ell_{m+1-i}.
$$
There is now only one non-degenerate element of colength $1$,
namely $w_ms_{\alpha}$ with $\alpha=\epsilon_1+\epsilon_m$.
We get
$$
(p^m-1) \lambda_1 [\Vc_{w_m}]\deg(\pi_{w_{m}})= (p+1) [\Vc_{w_{m+1}}]
\deg(\pi_{w_{m+1}}).
$$
Again by Lemma \ref{DegreelemmaD} we have
$\deg(\pi_{w_{m}})=\deg(\pi_{w_{m+2}})=1$.

Now take $w=w_{m+j}=s_{m-j}\ldots s_1$ with $j\geq 2$. The element
$v:=w_{\emptyset}w$ has an odd orbit of reduced orbit length $m-j+1$. We thus get
$$
(1-pv) \lambda_1= \ell_1+\sum_{i=1}^{m-j} p^i \epsilon_{m+2-j-i}.
$$
There is again only one non-degenerate element $ws_{\alpha}$ with
$\alpha=\epsilon_1-\epsilon_{m+1-j}$. Therefore
$\langle \alpha^{\vee},\lambda\rangle= (1-p^{j-1})$. We find
$$
(p^{m+1-j}+1) \lambda_1 [\Vc_{w_{m+j}}]
\deg(\pi_{w_{m+j}})=(p^{j-1}-1) [\Vc_{w_{m+j}}] \deg(\pi_{w_{m+j+1}}).
$$
By Lemma \ref{DegreelemmaD} we have
$\deg(\pi_{w_{m+j+1}})=p^{2j-2}+\cdots+p^{j}+p^{j-2}+\cdots+1$. Using the
factorization
$p^{2j-2}+\cdots+p^{j}+p^{j-2}+\cdots+1=(p^{j}+1)(p^{j-2}+p^{j-2}+\cdots+1)$ and
iterating we get the formula. Again the polynomiality is most easily verified by
explicit computation.
\end{proof}
\end{theorem}
\begin{remark}
This time formula iii) can be rewritten
\begin{displaymath}
\frac{1}{2}
\frac{p^{k-1}-1}{p^m+1}\left(\prod_{i=1}^{m-k+1}(p^i-1) \right)
\left[
\begin{array}{c}
m\\
k-1
\end{array}
\right]_{p^2}.
\end{displaymath}
\end{remark}
\begin{remark} It is not unreasonable to conjecture 
that the $\Vc_{w_k}$ are complete in the (open) moduli space.
for $k\geq 3$. Moreover, the class $\lambda_1$ is conjectured to be ample on the moduli space. In characteristic $0$ this follows from Baily-Borel.
[Recently this has been proved also in positive characteristic 
under mild conditions by Maulik \cite{M12} and Madapusi Pera 
\cite{MP12,MP13}.]
If this is true then the open strata $\sV_{w_{k}}$ for $k\geq 3$
are affine.
\end{remark}
\end{subsubsection}
\end{section}
\begin{section}{Applications}
We shall now discuss two applications both pertaining to the even case.
\begin{subsection}{(Quasi-)Elliptic fibrations with a section}

  If $X$ is a K3 surface and ${f}: {X}\to {\PP^1}$ is an elliptic (or possibly
  quasi-elliptic in characteristic $3$) fibration with a
  section $E\subset X$, then $E$ and a general fibre $F$ span a
  hyperbolic plane $\HH$ in $\NS(X)$ thus giving a $\HH$-marking of $X$. Let now
  $\sM^{es}$ be the stack of K3 surfaces together with a (quasi-)elliptic fibration
  (with base $\PP^1$) with a chosen section on it. As the choice of an ample
  line bundle is not part of the choices made let us take a moment to explain
  why this is an algebraic stack. We can consider the stack of K3-like surfaces
(i.e., a
  surface with rational double points only as singularities and whose minimal
  resolution is a K3 surface) with an (quasi-)elliptic fibration with a section and
  irreducible fibres which is smooth along the section. Three times a fibre plus
  the section is an ample divisor and hence the stack of such surfaces is
  algebraic. Then $\sM^{es}$ is the Artin-Brieskorn simultaneous resolution
  stack of it.
\begin{proposition}\label{sigma0neq10}
$\sM^{es}$ is of twisted even type of rank $20$ so that Theorem \ref{Twisted
  classes} applies with $m=10$. In particular $\sigma_0=10$ is not possible.
\begin{proof}
  We start by verifying the versality hypothesis for $\sM^{es}$. Let us
  therefore fix a geometric point $X\to{\rm Spec}({\bf k})$. Recall that deformations of
  K3 surfaces are unobstructed and the derivative of the period map
  $H^1(X,T_X)\to{\rm Hom}(H^0(X,\Omega^2_X),H^1(X,\Omega^1_X))$ is an isomorphism.
  Consider now the closed (formal) subscheme $A$ of some formal universal
  deformation ${\mathcal X} \to S$ of $X$ defined by the condition that the $\HH$-marking
  of $X$ extend over $A$.  Then $A$ is defined by two equations and its tangent
  space, as a subspace of $H^1(X,T_X)$ is given by the condition that $v\cdot
  c_1(\sL)=0\in H^2(X,\sO_X)$ for all line bundles $\sL\in \HH$ and where
  $c_1(\sL)\in H^1(X,\Omega^1_X)$ is the Hodge cohomology Chern class (induced
  by ${\mathrm{dlog}}:{\sO_X^*} \to {\Omega^1_X})$. As the degree of the polarization is
  prime to $p$, the class
 $c_1$ gives an injection $\HH\Tensor{\bf k}\hookrightarrow
  H^1(X,\Omega^1_X)$ and hence the codimension of $T_s(A)$ in $T_s(S)$, where
  $s$ is the closed point of $S$, is $2$ and hence $A$ is smooth. Furthermore,
  it also follows that $T_0(A)$ maps isomorphically onto
  ${\rm Hom}(H^0(X,\Omega^2_X),P)$, where $P$ is the primitive part of
  $H^1(X,\Omega^1_X)$. This gives the required versality for the stack of
  $\HH$-marked surfaces. What remains to show is that if the marking of $X$
  comes from a (quasi-)elliptic fibration with a section then so does any deformation
  of it. For the fibration we let $\sL=\sO_X(F)$ be the line bundle of a fibre
  $F$. Then $H^1(X,\sL)=0$ and hence for any extension of it (to some closed
  subscheme of $S$), its direct image will be a vector bundle $\sE$ of rank $2$
  which gives a map to the $\PP(\sE)$-bundle extending the (quasi-)elliptic fibration.
  Similarly, the section is a $(-2)$-curve $E$ and as $H^1(X,\sO_X(E))=0$ any
  extension of $\sO_X(E)$ will give an extension of the curve which then is a
  section. By construction both line bundles extend over $A$.

  As the discriminant of $\HH$ is $-1$ we get by Theorem \ref{Final/canonical
    filtrations type} that the Hodge discriminant of the primitive part is equal
  to $1$. Then from Proposition \ref{Odd Hodge discriminant} (and the fact that
  in the notations of that proposition $m=10$) we conclude that we are in
  the twisted case. Theorem \ref{Twisted classes} then gives the classes of the
  height and Artin invariant strata (together with the fact that $\sigma_0=10$
  is not possible).
\end{proof}
\end{proposition}
\begin{remark}
There is an alternative way of excluding $\sigma_0=10$ similar to the way Artin
excluded $\sigma_0=11$ for a general supersingular K3 surface. By
\cite{artin73::shafar+tate+k3} a supersingular (quasi-)elliptic K3 surface $X$ has
$\rho=22$ and by the fact that $\HH$ is unimodular we get that
$\mathrm{NS}(X)=\HH\perp P$. If $\sigma_0(X)=10$, then the scalar product on $P$
is divisible by $p$ and $P(1/p)$ is a unimodular even negative definite form of rank
$20$ which is not possible as its index, $20$, is not divisible by $8$. This
argument has the advantage of working also for $p=2$.

Since the first version of this paper was written two alternative proofs of 
Proposition \ref{sigma0neq10} have appeared, namely in the paper by Kondo and Shimada
\cite{KS}
and also in Liedtke's paper \cite{L13}.
\end{remark}
\end{subsection}
\begin{subsection}{The canonical double cover of an Enriques surface}
 We let $N$ be the lattice $E_{10}(-1)=\HH\perp E_8(-1)$ and we fix a chamber
  (inside of the positive cone) with respect to the roots of $N$ (see
  \cite[II:\S5]{cossec89::enriq+i} for a discussion of chambers in
  $E_{10}(-1)$). Let $\sM^E$ be the moduli stack of marked Enriques surfaces
  where a marking is an isometry between the standard Enriques lattice
  $N_{10}=\HH\perp E_8(-1)$ and the N\'eron-Severi group taking the fixed chamber
  into the ample cone of the N\'eron-Severi group. We can then construct
  $\sM^{E,d}\to \sM^{E}$, the moduli stack of canonical double covers of marked
  Enriques surfaces (i.e., while $\sM^E(S)$, for a scheme $S$, is the groupoid
  of families of marked Enriques surfaces over $S$, $\sM^{E,d}(S)$ is the
  groupoid of families of marked Enriques surfaces together with an unramified
  double cover of the Enriques surface which is fibrewise non-trivial).
\begin{remark}
  Note that $\sM^{E,d}\to \sM^{E}$ is not an isomorphism but rather a
  (non-trivial) $\ZZ/2$-gerbe. The non-triviality is reflected in the fact that
  given a family $X\to S$ of Enriques surfaces a canonical double cover is a
  double cover $X'\to X$ which is non-trivial over every geometric fibre over
  $S$. There is an obstruction in $H^2(S,\ZZ/2)$ which in general is non-zero to
  the existence of such a cover, making ``canonical double cover'' something of
  a misnomer.
\end{remark}
Pulling back the N\'eron-Severi group along the universal double cover 
${\mathcal X}'\to {\mathcal X}$ over
$\sM^{E,d}$ we get a marking by $N(2)$ of the family ${\mathcal X}'\to
\sM^{E,d}$ of K3 surfaces.
\begin{proposition}
$\sM^{E,d}$ is of twisted even type of rank $12$ so that Theorem \ref{Twisted
  classes} applies with $m=6$. In particular $\sigma_0=6$ is not possible.
\begin{proof}
  Again we start by verifying that $\sM^{E,d}$ fulfils the versality condition.
  We use the fact that $\sM^{E,d}$ also can be described as the stack of
  K3 surfaces together with $\iota$, a fixed point free involution. Its tangent
  space is then the space of linear maps $H^0(X,\Omega^1_X)\to H^1(X,\Omega_X)$
  commuting with the involution. Now, $\iota$ acts by $-1$ on
  $H^0(X,\Omega^1_X)$ and by $+1$ on $N(2)\Tensor{\bf k}$ and $-1$ on $P$ under the
  decomposition $H^1(X,\Omega_X)=N(2)\Tensor{\bf k}\perp P$ which gives what we want.

  The marking has discriminant $-2^{10}$ which is $-1$ up to squares just as in
  the previous example. Hence we are in the twisted even case (with $m=6$) and
  again Theorem \ref{Twisted classes} applies.
\end{proof}
\end{proposition}
\begin{remark}
Also in this case there is an arithmetic proof of the impossibility of
$\sigma_0=6$ (using as above that the scalar product on $P$ is divisible by $p$
and that its rank is not divisible by $8$). 
The proof does not extend to characteristic $2$ however as the
polarization is not of degree prime to $2$ (and the situation is in fact quite
different in characteristic $2$).
\end{remark}
\end{subsection}
\end{section}

\end{document}